\documentclass[reqno]{amsart}

\usepackage{hyperref}

\usepackage{amsmath}
\usepackage{amsfonts}
\usepackage{amsthm}
\usepackage{graphicx,epsfig}

\usepackage{color}

\numberwithin{equation}{section}

\newcommand{\supp}{\operatorname{supp}}
\newcommand{\Spec}{\operatorname{Spec}}

\newcommand{\Div}{{\operatorname{div}\,}}

\newcommand{\dist}{{\operatorname{dist}}}

\newcommand{\curl}{{\operatorname{curl}\,}}

\newcommand{\meas}{{\operatorname{meas}\,}}

\newtheorem{theorem}{Theorem}[section]
\newtheorem{thm}[theorem]{Theorem}
\newtheorem{lemma}[theorem]{Lemma}
\newtheorem{proposition}[theorem]{Proposition}
\newtheorem{prop}[theorem]{Proposition}
\newtheorem{definition}[theorem]{Definition}
\newtheorem{remark}[theorem]{Remark}

\newtheorem{assumption}[theorem]{Assumption}

\def \R {{\mathbb R}}

\title{Superconductivity in domains with corners}

\author{V.~Bonnaillie-No\"{e}l}
\author{S. Fournais}

\address[V.~Bonnaillie-No\"{e}l]{IRMAR, ENS Cachan Bretagne, Univ. Rennes 1, CNRS, UEB\\
av Robert Schuman, F-35170 Bruz, France}
\email{Virginie.Bonnaillie@Bretagne.ens-cachan.fr}

\address[S. Fournais]{Department of Mathematical Sciences, University of Aarhus, Ny Munke\-gade,
Building 1530,
DK-8000 Aarhus C, Denmark
}
\email{fournais@imf.au.dk}

\date{\today}

\begin{document}

\bibliographystyle{plain}

\begin{abstract}
We study the two-dimensional Ginzburg-Landau functional in a domain with corners for exterior magnetic field strengths near the critical field where the transition from the superconducting to the normal state occurs. We discuss and clarify the definition of this field and obtain a complete asymptotic expansion for it in the large $\kappa$ regime. Furthermore, we discuss nucleation of superconductivity at the boundary.
\end{abstract} 

\maketitle
\tableofcontents

\section{Introduction}
It is a well-known phenomenon that superconductors of Type II lose their superconducting properties when submitted to sufficiently strong external fields. The value of the external field where this transition takes place is usually called $H_{C_3}$, and is calculated as a function of a material-dependent parameter $\kappa$. The calculation of this critical field, $H_{C_3}$, for large values of $\kappa$ has been the focus of much activity \cite{BeSt}, \cite{LuPa1,LuPa2,LuPa3}, \cite{PiFeSt}, \cite{He-Mo} 
and \cite{He-Pan}. In the recent works \cite{FournaisHelffer3,FournaisHelffer4} the definition of 
$H_{C_3}$ in the case of samples of smooth cross section was clarified and it was realized that 
the critical field is determined completely by a linear eigenvalue problem.
The linear spectral problem has been studied in depth in the case of corners in \cite{Bo1,Bo2,BoDa}.
The objective of the present paper is to use the spectral information from \cite{BoDa} to carry through an  analysis similar to the one in \cite{FournaisHelffer3} in the case of corners.
Thereby we will in particular obtain: 1) A complete asymptotics of $H_{C_3}$ for large values of $\kappa$ in terms of linear spectral data, 2) Precise estimates on the location of nucleation of superconductivity for magnetic field strengths just below the critical field.\\
The case of corners of angle $\pi/2$ has been studied in \cite{Jad,Pancorner}. Our results are more precise---even for those angles---and we study more general domains.

We will work in the Ginzburg-Landau model.
Let $\Omega \subset {\mathbb R}^2$ be a bounded simply connected domain with Lipschitz boundary.
The Ginzburg-Landau functional is given by
\begin{multline}
\label{eq:GL_F}
{\mathcal E}[\psi,{\bf A}] = {\mathcal
E}_{\kappa,H}[\psi,{\bf A}]  =
\int_{\Omega} \Big\{ |p_{\kappa H {\bf A}}\psi|^2 
- \kappa^2|\psi|^2
+\frac{\kappa^2}{2}|\psi|^4\Big\} \,dx\\
+ \kappa^2 H^2 \int_{{\mathbb R}^2}
|\curl {\bf A} - 1|^2\,dx\;,
\end{multline}
with
$\psi \in W^{1,2}(\Omega;{\mathbb C})$, ${\bf A}$ in the space $\dot{H}^1_{{\bf F}, \Div}$ that we will define below, and where $p_{{\bf A}} = (-i\nabla- {\bf A})$. 
Notice that the second integral in \eqref{eq:GL_F} is over the entire space, ${\mathbb R}^2$, whereas the first integral is only over the domain $\Omega$.

Formally the functional is gauge invariant. In order to fix the gauge, we will impose that vector fields ${\bf A}$ have vanishing divergence. Therefore, a good choice for the variational space for ${\bf A}$ is 
\begin{align}
\dot{H}^1_{{\bf F}, \Div} = {\bf F} + \dot{H}^1_{\Div}, \end{align}
where
\begin{align*}
\dot{H}^1_{\Div} = \{ {\bf A} \in \dot{H}^1({\mathbb R}^2, {\mathbb R}^2) \,\big |\, \Div {\bf A} = 0 \}\;.
\end{align*}
Furthermore ${\bf F}$ is the vector potential giving constant magnetic field
\begin{align}
{\bf F}(x_1,x_2) = \tfrac{1}{2}(-x_2, x_1),
\end{align}
and we use the notation $\dot{H}^1({\mathbb R}^2)$ for the homogeneous Sobolev spaces, i.e. the closure of $C_0^{\infty}({\mathbb R}^2)$ under the norm
$$
f \mapsto \| f \|_{\dot{H}^1} = \| \nabla f \|_{L^2}.
$$
Any square integrable magnetic field $B(x)$ can be represented by a vector field ${\bf A} \in \dot{H}^1_{\Div}$.

Minimizers, $(\psi, {\bf A}) \in W^{1,2}(\Omega) \times \dot{H}^1_{{\bf F}, \Div}$, of the functional ${\mathcal E}$ have to satisfy the Euler-Lagrange equations:
\begin{subequations}
\label{eq:GL}
\begin{align}
\label{eq:equationA}
p_{\kappa H {\bf A}}^2\psi =
\kappa^2(1-|\psi|^2)\psi & \quad \text{ in } \quad\Omega\;,\\
\label{eq:equationB}
\curl^2 {\bf A} =\big\{-\tfrac{i}{2\kappa H}(\overline{\psi} \nabla
\psi - \psi \nabla \overline{\psi}) - |\psi|^2 {\bf A}\big\} 1_{\Omega}(x)
&\quad \text{ in } \quad {\mathbb R}^2 \;,\\
(p_{\kappa H {\bf A}} \psi) \cdot \nu = 0 &\quad \text{ on } \quad \partial\Omega\;.
\end{align}
\end{subequations}
It is standard to prove that for all $\kappa, H>0$, the functional ${\mathcal E}_{\kappa,H}$ has a minimizer. An important result by Giorgi and Phillips, \cite{Giorgi-Phillips}, states that for $\kappa$ fixed and $H$ sufficiently large (depending on $\kappa$), the unique solution of \eqref{eq:GL} (up to change of gauge) is the pair $(\psi, {\bf A}) = (0, {\bf F})$. Since $\psi$ is a measure of the superconducting properties of the state of the material and ${\bf A}$ is the corresponding configuration of the magnetic vector potential, the result of Giorgi and Phillips reflects the experimental fact that superconductivity is destroyed in a strong external magnetic field. 

We define the lower critical field $\underline{H}_{C_3}$ as the value of $H$ where this transition takes place:
\begin{align}
\label{eq:DefHC3}
\underline{H}_{C_3}(\kappa) = \inf\{ H>0 \;:\; (0, {\bf F}) \text{ is a minimizer of } {\mathcal E}_{\kappa,H}\}\;.
\end{align}
However, it is far from obvious from the functional that the transition takes place at a unique value of $H$---there could be a series of transitions back and forth before the material settles definitely for the normal state, $(0,{\bf F})$. Therefore, we introduce a corresponding upper critical field
\begin{align}
\label{eq:DefUpper}
\overline{H}_{C_3}(\kappa) &=  \inf\{ H>0 \;:\; (0, {\bf F}) \text{ is the unique minimizer of } {\mathcal E}_{\kappa,H'}  \text{ for all } H'>H\} \;.
\end{align}
Part of our first result, Theorem~\ref{thm:EqualHC3} below, is that the above definitions coincide for large $\kappa$.

Let us introduce some spectral problems. For $B\geq 0$ and a (sufficiently regular) domain $\Omega \subset {\mathbb R}^2$, we can define a quadratic form 
\begin{align}
\label{eq:Form}
Q[u] = Q_{\Omega,B}[u] = \int_{\Omega} |(-i\nabla - B {\bf F})u|^2\,dx,
\end{align}
with form domain $\{ u \in L^2(\Omega) \,\big |\, (-i\nabla - B {\bf F})u \in L^2(\Omega) \}$. The self-adjoint operator associated to this closed quadratic form will be denoted by ${\mathcal H}(B) = {\mathcal H}_{\Omega}(B)$. Notice that since the form domain is maximal, the operator ${\mathcal H}_{\Omega}(B)$ will correspond to Neumann boundary conditions.
We will denote the $n$'th eigenvalue of ${\mathcal H}(B)$ (counted with multiplicity) by
$\lambda_n(B)= \lambda_{n,\Omega}(B)$,
in particular,
$$
\lambda_1(B) = \lambda_{1,\Omega}(B):= \inf \Spec {\mathcal H}_{\Omega}(B).
$$

The case where $\Omega$ is an angular sector in the plane will provide important special models for us. Define, for $0<\alpha\leq 2\pi$,
$$
\Gamma_{\alpha} :=\{ z = r(\cos \theta, \sin \theta) \in {\mathbb R}^2 \,\big| \, r \in (0,\infty), |\theta| < \alpha/2 \}.
$$
Since this domain is scale invariant one easily proves that
$$
\Spec {\mathcal H}_{\Gamma_{\alpha}}(B) = B \Spec {\mathcal H}_{\Gamma_{\alpha}}(1).
$$
Therefore, we set $B=1$ and define
\begin{align}
\mu_1(\alpha) = \lambda_{1,\Gamma_{\alpha}}(B=1).
\end{align}
The special case of $\alpha = \pi$, i.e. the half plane, has been studied intensively. In compliance with standard notation, we therefore also write
$$
\Theta_0 := \mu_1(\alpha=\pi).
$$
It is known that the numerical value of $\Theta_0$ is $\Theta_0 = 0.59....$.
\begin{remark}~\\
It is believed---and numerical evidence exists (cf. \cite{AlBo,BDMV} and Figure~\ref{FigMuAlpha}) to support this claim---that $\alpha \mapsto \mu_1(\alpha)$ is a strictly increasing function on $[0,\pi]$ and constant equal to $\Theta_0$ on $[\pi, 2\pi]$. If this belief is proved, then the statement of our Assumption~\ref{ass:Angles} below can be made somewhat more elegantly.
\begin{figure}[h!]\begin{center}
\includegraphics[height=7cm]{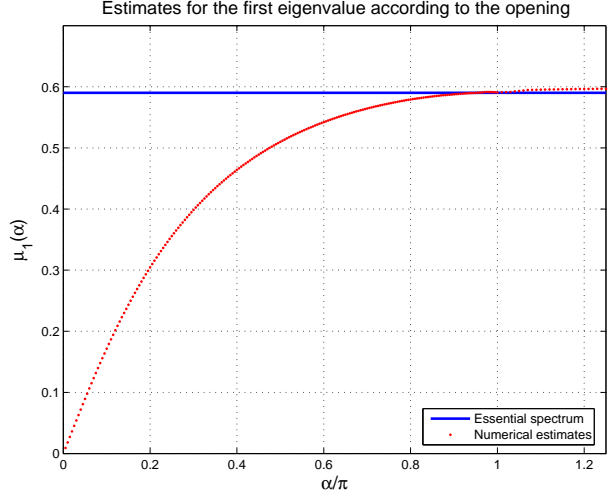}
\caption{$\mu_1(\alpha)$ vs. $\alpha/ \pi$ for $\alpha\in [0, 1.25 \pi]$.}\label{FigMuAlpha}
\end{center}\end{figure}
\end{remark}

We consider $\Omega$ a domain whose boundary is a curvilinear polygon in the sense given by Grisvard, see Definition~\ref{def.Gr}. 

\begin{definition}[cf. {\cite[p.34--42]{Gr}}] \label{def.Gr}~\\
Let $\Omega$ be a bounded open subset of $\mathbb R^2$. We say that the boundary $\Gamma$ is a
(smooth) curvilinear polygon, if for every $x\in\Gamma$ there exists a neighborhood  $V$ of $x$ in $\mathbb R^2$ and a mapping $\psi$ from $V$ in $\mathbb R^2$ such that
\begin{enumerate}
\item $\psi$ is injective,
\item $\psi$ together with $\psi^{-1}$ (defined on $\psi(V)$) belongs to the class $C^{\infty}$,
\item $\Omega\cap V$ is either
$\{y\in\Omega\,|\, \psi_{2}(y)<0\}$, $\{y\in\Omega\,|\, \psi_{1}(y)<0\mbox{ and }\psi_{2}(y)<0\}$,
or $\{y\in\Omega\,|\, \psi_{1}(y)<0\mbox{ or }\psi_{2}(y)<0\}$, 
where $\psi_{j}$ denotes the components of $\psi$.
\end{enumerate}
\end{definition}

From now on, we consider a bounded open subset $\Omega\subset\mathbb R^2$, whose boundary is a curvilinear polygon of class $C^{\infty}$. 
The boundary of such a domain will be a piecewise smooth curve $\Gamma$.
We denote the (minimal family of) smooth curves which make up the boundary by $\overline \Gamma_{j}$ for $j=1,\ldots,N$. The curve $\overline \Gamma_{j+1}$ follows $\overline\Gamma_{j}$ according to a positive orientation, on each connected component of $\Gamma$. We denote by ${\mathsf{s}}_{j}$ the vertex which is the end point of $\overline \Gamma_{j}$. We define a vector field $\nu_{j}$ on a neighborhood of $\overline\Omega$, which is the unit normal a.e. on $\Gamma_{j}$. 

We will work under the following assumption on the domain.

\begin{assumption}\label{ass:Angles}~\\
The domain $\Omega$ has curvilinear polygon boundary and denote the set of vertices by $\Sigma$. 
We suppose that $N:= |\Sigma| \neq 0$.
We denote by $\alpha_{\mathsf{s}}$ the angle at the vortex $\mathsf{s}$ (measured towards the interior). 
We suppose that $\mu_1(\alpha_{\mathsf{s}})< \Theta_0$ for all $\mathsf{s}\in\Sigma$,
and define $\Lambda_1 := \min_{\mathsf{s} \in \Sigma} \mu_1(\alpha_{\mathsf{s}})$. 
We also assume that $\alpha_{\mathsf{s}} \in (0, \pi)$ for all $\mathsf{s}\in\Sigma$.
\end{assumption}

Under this assumption we resolve the ambiguity of definition of $H_{C_3}(\kappa)$ and derive a complete asymptotics in terms of spectral data.

\begin{thm}\label{thm:EqualHC3}~\\
Suppose that $\Omega$ is a bounded, simply-connected domain satisfying Assumption~\ref{ass:Angles}.
Then there exists $\kappa_0>0$ such that if $\kappa \geq \kappa_0$ then the equation
$$
\lambda_{1,\Omega}(\kappa H) = \kappa^2,
$$
has a unique solution $H=H_{C_3}^{\rm lin}(\kappa)$. 
Furthermore, if $\kappa_0$ is chosen sufficiently large, then for $\kappa \geq \kappa_0$, the critical fields defined in \eqref{eq:DefHC3}, \eqref{eq:DefUpper} coincide and satisfy
\begin{align}
\label{eq:EqualFields}
{\underline H}_{C_3}(\kappa) = {\overline H}_{C_3}(\kappa) = H_{C_3}^{\rm lin}(\kappa).
\end{align}
Finally, the critical field has a complete asymptotic expansion in powers of $\kappa^{-1}$: There exists $\{\eta_j\}_{j=1}^{\infty} \subset {\mathbb R}$ such that
\begin{align}
\label{eq:Asymptotics}
{\underline H}_{C_3}(\kappa) = \frac{\kappa}{\Lambda_{1}}\Big(1 + \sum_{j=1}^{\infty} \eta_j \kappa^{-j}\Big),\quad \text{ for } \quad \kappa \rightarrow \infty,
\end{align}
in the sense of asymptotic series.
\end{thm}

\begin{remark}~\\
The result analogous to Theorem~\ref{thm:EqualHC3} for smooth domains (i.e. for $\Sigma = \emptyset$) has been established in \cite{FournaisHelffer3,FournaisHelffer4}. Notice however that the form of the asymptotics \eqref{eq:Asymptotics} depends on the existence of a vortex and is more complicated in the case of smooth domains.
\end{remark}

Once Theorem~\ref{thm:EqualHC3} is established it makes sense, for large values of $\kappa$, to talk of {\it the} critical field that we will denote by $H_{C_3}(\kappa)$  ($= {\underline H}_{C_3}(\kappa) = {\overline H}_{C_3}(\kappa)$).

In the case of regular domains (without corners) one has the asymptotics
(see \cite{LuPa1}, \cite{PiFeSt}, \cite{He-Mo} and \cite{He-Pan}),
$$
H_{C_3}(\kappa) = \frac{\kappa}{\Theta_0} + {\mathcal O}(1),
$$
where the leading correction depends on the maximal curvature of the boundary.
We observe that the corners---which can be seen as points where the curvature is infinite---change the leading order term of $H_{C_3}(\kappa)$.
Thus there is a large parameter regime of magnetic field strengths, ${\kappa}/{\Theta_0} \ll H \leq H_{C_3}(\kappa)$, where superconductivity in the sample must be dominated by the corners.
Our next two results make this statement precise.
First we prove Agmon type estimates, for the minimizers of the non-linear Ginzburg-Landau functional, which describe how superconductivity can nucleate successively in the corners, ordered according to  their spectral parameter $\mu_1(\alpha_{\mathsf s})$.

\begin{thm}\label{thm:AgmonCorners}~\\
Suppose that $\Omega$ satisfies Assumption~\ref{ass:Angles},
let $\mu>0$ satisfy $\min_{\mathsf{s}\in\Sigma}\mu_1(\alpha_{\mathsf{s}}) < \mu < \Theta_0$ and define
$$
\Sigma' := \{ {\mathsf s} \in \Sigma\,\big|\, \mu_1(\alpha_{\mathsf s}) \leq \mu \}.
$$
There exist constants $\kappa_0, M, C, \epsilon>0$ such that if 
$$
\kappa \geq \kappa_0, \quad\quad\quad
\frac{H}{\kappa} \geq \mu^{-1},
$$
and $(\psi, {\bf A})$ is a minimizer of ${\mathcal E}_{\kappa,H}$, then
\begin{align*}
\int_{\Omega} e^{\epsilon \sqrt{\kappa H} \dist(x, \Sigma')}
\bigg( |\psi(x)|^2 &+ \frac{1}{\kappa H} |p_{\kappa H {\bf A}} \psi(x)|^2\bigg)\,dx \\
&\leq C \int_{\{x: \sqrt{\kappa H}\dist(x,\Sigma') \leq M\}} |\psi(x)|^2\,dx.
\end{align*}
\end{thm}

Finally we discuss leading order energy asymptotics in the parameter regime dominated by the corners, i.e. ${\kappa}/{\Theta_0} \ll H \leq {\underline H}_{C_3}(\kappa).$
The result below, Theorem~\ref{thm:EnergyCorners}, can be seen as a partial converse to Theorem~\ref{thm:AgmonCorners} in that all corners which are spectrally permitted will contribute to the leading order of the ground state energy.

One can imagine an interaction between corners with the same spectral parameter, i.e. with the same angle $\alpha$. This would be a tunnelling type effect and has much lower order. We refrain from a detailed study of such an interaction, since that would be far out of the scope of the present paper.

The ground state energy will be given to leading order by decoupled model problems in angular sectors. It may be slightly surprising to notice that these model problems remain {\it non-linear}.

Let $\alpha \in(0,\pi)$ be such that $\mu_1(\alpha) < \Theta_0$. (Remember that it follows from \cite{Bo2} that $\mu_1(\alpha)<\Theta_0$ for $\alpha\in (0,\frac\pi 2]$ and that numerical evidence suggests this to be the case in the entire interval $\alpha\in (0,\pi )$.)

Define, for $\mu_1, \mu_2 > 0$, the following functional $J_{\mu_1, \mu_2}^{\alpha}$,
\begin{align}
\label{eq:Jmu}
J_{\mu_1, \mu_2}^{\alpha}[\psi] =
\int_{\Gamma_{\alpha}} \left\{ |(-i\nabla -{\bf F})\psi|^2 - \mu_1 |\psi |^2 + \frac{\mu_2}{2}|\psi|^4\right\}\,dx,
\end{align}
with domain 
$\{ \psi \in L^2(\Gamma_{\alpha})\,| \,(-i\nabla -{\bf F})\psi \in L^2(\Gamma_{\alpha})\}$.
Define also the corresponding ground state energy
$$
E_{\mu_1, \mu_2}^{\alpha}:= \inf J_{\mu_1, \mu_2}^{\alpha}[\psi] .
$$
The main result on the ground state energy of the Ginzburg-Landau functional in the parameter regime dominated by the corners is the following.

\begin{thm}\label{thm:EnergyCorners}~\\
Suppose $\frac{\kappa}{H(\kappa)} \rightarrow \mu \in {\mathbb R}_{+}$ as $\kappa \rightarrow \infty$, where $\mu < \Theta_0$.
Let $(\psi, {\bf A})= (\psi, {\bf A})_{\kappa,H(\kappa)}$ be a minimizer of ${\mathcal E}_{\kappa, H(\kappa)}$.

Then
\begin{align}
\label{eq:EnergyAsymp}
{\mathcal E}_{\kappa, H(\kappa)}[\psi, {\bf A}] \rightarrow \sum_{\mathsf{s}\in\Sigma} E^{\alpha_{\mathsf{s}}}_{\mu,\mu},
\end{align}
as $\kappa \rightarrow \infty$.
\end{thm}

\begin{remark}~\\
Proposition~\ref{prop:AngularNL} below states that
$E^{\alpha_{\mathsf{s}}}_{\mu,\mu}=0$ unless $\mu_1(\alpha_{\mathsf{s}}) < \mu$, so only corners satisfying this spectral condition contribute to the ground state energy in agreement with the localization estimate from Theorem~\ref{thm:AgmonCorners}.
\end{remark}

\section{Spectral analysis of the linear problem}
\subsection{Monotonicity of $\lambda_1(B)$}~\\
In this subsection we will prove that $B \mapsto \lambda_1(B)$ is increasing for large $B$. Thereby we will have proved the first statement of Theorem~\ref{thm:EqualHC3} (see Propositions~\ref{prop.lambda1'} and \ref{prop.lambda1} below). Furthermore, Lemma~\ref{lem:Asymptotics} establishes the form of the asymptotics of $H^{\rm lin}_{C_3}(\kappa)$.

In \cite{BoDa} the asymptotics of $\lambda_1(B)$ was effectively calculated to any order. Let us recall their results. 
\begin{definition}\label{def.polyg}~\\
Let $\Omega$ be a bounded curvilinear polygon. We denote by 
\begin{itemize}
\item $\Lambda_{n}$ the $n$-th eigenvalue of the model operator $\oplus_{\mathsf{s}\in\Sigma}Q^{\alpha_{\mathsf{s}}}$ where $Q^{\alpha_{\mathsf{s}}}$ is the magnetic Neumann Laplacian $(-i\nabla -{\bf F})^2$ on the infinite angular sector of opening $\alpha_{\mathsf{s}}$.
In particular, $\Lambda_1 = \min_{\mathsf{s}\in\Sigma}\mu_1(\alpha_{\mathsf{s}})$,
\item $K_{\Omega}$ the largest integer $K$ such that $\Lambda_{K}<\Theta_{0}$,
\item $\mu^{(n)}(h)$ the $n$-th smallest eigenvalue of the magnetic Neumann Laplacian $(-ih\nabla -{\bf F})^2$ on $\Omega$.
\end{itemize}
\end{definition}
\begin{theorem}[\cite{BoDa} Theorem 7.1]\label{th.polyg}~\\
Let $n\leq K_{\Omega}$. There exists $h_{0}>0$ and $(m_{j})_{j\geq 1}$ such that for any $N>0$ and $h\leq h_0$,
$$\mu^{(n)}(h)=h\Lambda_{n}+h\sum_{j=1}^Nm_{j}h^{j/2}+\mathcal{O}(h^{\frac{N+1}{2}}).$$
Furthermore, if $\Omega$ is a bounded convex polygon (i.e.~has straight edges), then for any $n\leq K_{\Omega}$, there exists $r_{n}>0$ and for any $\varepsilon>0$, $C_{\varepsilon}>0$ such that
$$\left|\mu^{(n)}(h)-h\Lambda_{n}\right|\leq C_\varepsilon \exp\left(-\frac{1}{\sqrt h}(r_{n}\sqrt{\Theta_0-\Lambda_{n}-\varepsilon})\right).$$ 
\end{theorem}
Recall the notation ${\mathcal H}(B)$, $\lambda_n(B)$ introduced after \eqref{eq:Form}.
By a simple scaling, we get
\begin{align}
\label{eq:scalingh-B}
\lambda_{n}(B)=B^2\mu^{(n)}(B^{-1}),\quad\forall n.
\end{align}
Let us make more precise the behavior of $\lambda_{1}(B)$ as $B$ is large. For this, we define the left and right derivatives of $\lambda_{1}(B)$:
\begin{equation}
\lambda_{1,\pm}'(B):=\lim_{\varepsilon\to 0^\pm}\frac{\lambda_{1}(B+\varepsilon)-\lambda_{1}(B)}{\varepsilon}.
\end{equation}
\begin{proposition}\label{prop.lambda1'}~\\
The limits of $\lambda_{1,+}'(B)$ and $\lambda_{1,-}'(B)$ as $B\to +\infty$ exist, are equal and we have
$$\lim_{B\to+\infty}\lambda_{1,+}'(B)=\lim_{B\to+\infty}\lambda_{1,-}'(B)=\Lambda_{1}.$$
Therefore, $B\mapsto\lambda_{1}(B)$ is strictly increasing for large $B$.
\end{proposition}
\begin{proof}~\\
This proof is similar to that of \cite{FournaisHelffer3}.\\
Let $B\geq 0$ and let $n$ be the degeneracy of $\lambda_{1}(B)$. There exist $\varepsilon>0$, $2n$ analytic functions $\phi_{j}$ and $E_{j}$, $j=1,\ldots,n$ defined from $(B-\varepsilon,B+\varepsilon)$ into 
$H^2(\Omega)\setminus\{0\}$ and $\mathbb{R}$ respectively, such that
$$
\mathcal{H}(\beta)\phi_{j}(\beta)=E_{j}(\beta)\phi_{j}(\beta),\quad E_{j}(B)=\lambda_{1}(B),
$$
and such that $\{\phi_j(B)\}$ are linearly independent. 
If $\varepsilon$ is small enough, there exist $j_{+}$ and $j_{-}$ in $\{1,\ldots,n\}$ such that
\begin{eqnarray*}
\mbox{for }\beta\in(B,B+\varepsilon),&& E_{j,+}(\beta)=\min_{j\in\{1,\ldots,n\}}E_{j}(\beta),\\
\mbox{for }\beta\in(B-\varepsilon,B),&& E_{j,-}(\beta)=\min_{j\in\{1,\ldots,n\}}E_{j}(\beta).
\end{eqnarray*}
By first order perturbation theory, the derivatives $\lambda_{1,\pm}'(B)$ can be rewritten
$$\lambda_{1,\pm}'(B)=-2\Re\langle \phi_{j_{\pm}}(B),{\bf F}\cdot(-i\nabla-B{\bf F})\phi_{j_{\pm}}(B)\rangle.$$
We deduce, for $\varepsilon>0$,
\begin{eqnarray*}
\lambda_{1,+}'(B)
&=&\frac1\varepsilon \langle \phi_{j_{+}}(B),(\mathcal{H}(B+\varepsilon)-\mathcal{H}(B)-\varepsilon^2{\bf F}^2)\phi_{j_{+}}(B)\rangle\\
&\geq &\frac{\lambda_{1}(B+\varepsilon)-\lambda_{1}(B)}\varepsilon-\varepsilon
\|{\bf F}\|^2_{L^\infty(\Omega)}.
\end{eqnarray*}
Using Theorem~\ref{th.polyg}, we deduce that 
$$
\lambda_{1,+}'(B)\geq \Lambda_{1} +
m_1 \frac{\sqrt{B+\varepsilon} -\sqrt{B}}{\varepsilon} + \varepsilon^{-1} {\mathcal O}(B^{-1/2})
-\varepsilon \|{\bf F} \|^2_{L^\infty(\Omega)}.
$$
Thus,
$$\liminf_{B\to\infty}\lambda_{1,+}'(B)\geq \Lambda_{1}-\varepsilon \|{\bf F}\|^2_{L^\infty(\Omega)}.$$
Since $\varepsilon$ is arbitrary, we have $$\liminf_{B\to\infty}\lambda_{1,+}'(B)\geq \Lambda_{1}.$$
Taking $\varepsilon<0$, we obtain by a similar argument, 
$$\liminf_{B\to\infty}\lambda_{1,-}'(B)\leq \Lambda_{1}.$$
The two last inequalities and the relation $\lambda_{1,+}'(B)\leq\lambda_{1,-}'(B)$ achieve the proof.
\end{proof}

We are now able to prove the following proposition.
\begin{proposition}\label{prop.lambda1}~\\
The equation in $H$ 
$$\lambda_{1}(\kappa H)=\kappa^2$$ 
has a unique solution $H(\kappa)$ for $\kappa$ large enough.
\end{proposition}
\begin{proof}~\\
According to Proposition~\ref{prop.lambda1'}, there exists $B_0>0$ such that $\lambda_{1}$ is a strictly increasing continuous function from $[B_{0},+\infty)$ onto $[\lambda_{1}(B_{0}),+\infty)$. 
By choosing $B_0$ sufficiently large, we may assume that $\lambda_1(B) < \lambda_{1}(B_{0})$ for all $B<B_0$.
Let $\kappa_{0}=\sqrt{\lambda_{1}(B_{0})}$, then, for any $B>B_{0}$, the equation
$$\lambda_{1}(\kappa H)=\kappa^2$$
has a unique solution $H=\lambda_{1}^{-1}(\kappa^2)/\kappa$ with $\lambda_{1}^{-1}$ the inverse function of $\lambda_{1}$ defined on $[\lambda_{1}(B_{0}),+\infty)$.
\end{proof}

\begin{lemma}\label{lem:Asymptotics}~\\
Let $H=H^{\rm lin}_{C_3}(\kappa)$ be the solution to the equation
$$\lambda_{1}(\kappa H)=\kappa^2$$
given by Proposition~\ref{prop.lambda1}.
Then there exists a real valued sequence $(\eta_{j})_{j\geq 1}$ such that
\begin{equation}\label{eq.Hkappa}
H^{\rm lin}_{C_3}(\kappa)=\frac{\kappa}{\Lambda_{1}}\Big(1+\sum_{j=1}^\infty \eta_{j}\kappa^{-j}\Big),
\end{equation}
(in the sense of asymptotic series) 
with $\Lambda_{1}=\min_{\mathsf{s}\in\Sigma}\mu_{1}(\alpha_{\mathsf{s}})$ introduced in Definition~\ref{def.polyg}.
\end{lemma}
\begin{proof}~\\
By Theorem~\ref{th.polyg} and \eqref{eq:scalingh-B} there exists a sequence $(m_{k})_{k\geq 1}$ such that, for any $N\in\mathbb{N}$,
\begin{align}
\label{eq:EigenvalueAsymp}
\lambda_{1}(B)=\Lambda_{1} B+B\sum_{k=1}^N m_{k} B^{-k/2}+\mathcal{O}(B^{\frac{-N+1}{2}})\quad\mbox{ as }B\to +\infty.
\end{align}
We compute with the {\it Ansatz} for $H(\kappa)$ given by \eqref{eq.Hkappa}~:
\begin{eqnarray*}
\lambda_{1}(\kappa H)&\sim&
\Lambda_{1}\kappa H+\kappa H\sum_{k\geq 1} m_{k} (\kappa H)^{-k/2}\\
&\sim& \kappa^2 \big(1+\sum_{j=1}^\infty \eta_{j}\kappa^{-j}\big)
+\sum_{k\geq 1}m_{k}\frac{\kappa^{2-k}}{ \Lambda_{1}^{1-k/2} } 
\big(1+\sum_{j=1}^\infty \eta_{j}\kappa^{-j}\big)^{1-k/2}\\
&=&\kappa^2 + \big(\eta_{1}+\frac{m_{1}}{\sqrt{\Lambda_{1}}}\big)\kappa+
\big(\eta_{2}+\frac{m_{1}}{\sqrt{\Lambda_{1} }}\frac{\eta_{1}}{2}+m_{2}\big)+
\ldots\\
&=&\kappa^2+\kappa^2\sum_{j\geq 1}(\eta_{j}+\tilde m_{j})\kappa^{-j},
\end{eqnarray*}
where the coefficients $\tilde m_{j}$ only depend on the $\eta_{k}$ for $k<j$. Thus, the form \eqref{eq.Hkappa} admits a solution in the sense of asymptotic series.
It is an easy exercise to prove that $H^{\rm lin}_{C_3}(\kappa)$ is equivalent to this series.
\end{proof}

\subsection{Agmon estimates near corners for the linear problem}~\\
If $\phi \in C_0^{\infty}(\Omega)$ (i.e. with support away from $\partial \Omega$) it is a simple calculation to prove that
\begin{align}
\int_{\Omega} |(-i\nabla -{\bf A})\phi|^2\,dx \geq \int_{\Omega} \curl {\bf A} |\phi|^2\,dx.
\end{align}
In particular, for ${\bf A} = B{\bf F}$,
\begin{align}
\label{eq:LowerInterior}
Q_{\Omega,B}[\phi] \geq B \|\phi\|^2.
\end{align}
Using the technique of Agmon estimates (\cite{Ag,He2}) one can combine the upper and lower bounds
\eqref{eq:EigenvalueAsymp} and \eqref{eq:LowerInterior} to obtain exponential localization near the boundary for ground state eigenfunctions of ${\mathcal H}(B)$. 
For completeness we give the following theorem (without proof---we will give the proof of similar {\it non-linear} estimates below), though we will not need the result here.
\begin{thm}~\\
Let $\psi_B$ be the ground state eigenfunction of ${\mathcal H}(B)$.
Then there exist constants $\epsilon, C, B_0>0$ such that
$$
\int e^{\epsilon\sqrt{B}\dist(x,\partial\Omega)}\big\{ |\psi_B(x)|^2 + B^{-1} |p_{B{\bf F}}\psi_B(x)|^2\big\}\,dx
\leq C \| \psi_B\|_2^2\;,
$$
for all $B\geq B_0$.
\end{thm}
In order to prove exponential localization near the corners for minimizers of ${\mathcal E}_{\kappa,H}$ we will need the operator inequality \eqref{eq:LowerPotential} below (compare to \eqref{eq:LowerInterior}).
\begin{thm}
\label{thm:LowerBoundPotentialCorners}~\\
Let $\delta>0$. Then there exist constants $M_0, B_0>0$ such that if $B \geq B_0$ then ${\mathcal H}(B)$ satisfies the operator inequality
\begin{align}
\label{eq:LowerPotential}
{\mathcal H}(B) \geq U_B,
\end{align}
where $U_B$ is the potential given by
$$
U_B(x):= \begin{cases}
(\mu_1(\alpha_{\mathsf{s}}) - \delta)B,& \dist(x, {\mathsf{s}}) \leq M_0/\sqrt{B},\\
(\Theta_0 - \delta)B, & \dist(x, \Sigma) > M_0/\sqrt{B},\; \dist(x, \partial \Omega) \leq M_0/\sqrt{B},\\
(1-\delta)B,& \dist(x, \partial \Omega) > M_0/\sqrt{B}.
\end{cases}
$$
\end{thm}

\begin{proof}[Proof of Theorem~\ref{thm:LowerBoundPotentialCorners}]~\\
Let $\chi_1 \in C^{\infty}({\mathbb R})$ be non-increasing and satisfy $\chi_1(t) = 1$ for $t\leq 1$, $\chi_1(t)=0$ for $t\geq 2$. 

Define, for $L,M,B>0$,
\begin{align*}
\chi^{\rm cor}_M(x) &:= \chi_1(\sqrt{B} \dist(x, \Sigma)/M),\\
\chi^{\rm bd}_M(x) &:= \sqrt{(1-\chi_1^2)\big(\sqrt{B}  \dist(x, \Sigma)/M\big)} \times
\chi_1\big(\sqrt{B} L \dist(x, \partial \Omega)/(2M)\big)\\
\chi^{\rm int}_M(x) &:= \sqrt{(1-\chi_1^2)\big(\sqrt{B} L \dist(x, \partial \Omega)/(2M)\big)} .
\end{align*}
The parameter $L$ will be fixed. It is chosen sufficiently large that $\supp \chi^{\rm bd}_M$ consists of $N$ (the number of smooth boundary curves) disjoint components (lying along each smooth boundary piece) when $\sqrt{B}/M$ is large.

Using the IMS-formula we can write for any $\phi \in H^1(\Omega)$,
\begin{align}
\label{eq:IMS}
Q_{\Omega,B}[\phi] \geq
Q_{\Omega,B}[\chi^{\rm cor}_{M} \phi] + Q_{\Omega,B}[\chi^{\rm bd}_{M} \phi] 
+ Q_{\Omega,B}[\chi^{\rm int}_{M} \phi] 
-C \frac{B}{M^2} \| \phi\|^2,
\end{align}
for some constant $C>0$ independent of $M$, $B$ and $\phi$.

We will estimate each term of \eqref{eq:IMS} by using successively results for the first eigenvalue of the Schr\"odinger operator in a domain with one corner, in a smooth domain and in the entire plane.

Since $\chi^{\rm int}_M\phi$ has compact support in $\Omega$, we get (see \eqref{eq:LowerInterior})
\begin{align}
\label{eq:INT}
Q_{\Omega,B}[\chi^{\rm int}_M\phi]
=Q_{\mathbb{R}^2,B}[\chi^{\rm int}_M\phi] \geq B \| \chi^{\rm int}_M\phi \|^2.
\end{align}

For the corner contribution  and boundary contribution, we will use the estimates in angular sectors and regular domains obtained in \cite{Bo2,He-Mo}.\\
For any corner $\mathsf{s}\in\Sigma$, we define a domain $\Omega_{\mathsf{s}}$ such that $\Omega\cap B(\mathsf{s},\varepsilon) =  \Omega_{\mathsf{s}}\cap B(\mathsf{s},\varepsilon)$ for $\varepsilon$ small enough ($\varepsilon< {\rm dist}(\mathsf{s},\Sigma\setminus\{\mathsf{s}\})$) and its boundary is $C^\infty$ except in $\mathsf{s}$. 
Let $\mathsf{s}_{-}$ and $\mathsf{s}_{+}$ be the neighbor vertices of $\mathsf{s}$ (if they exist). 
We define two regular domains $\Omega_{\mathsf{s}}^-$ and $\Omega_{\mathsf{s}}^+$ such that there exists $\varepsilon>0$ with $\Omega\cap B(x,\varepsilon) = \Omega_{\mathsf{s}}^\pm \cap B(x,\varepsilon)$ for any $x\in \{y\in \Gamma_{\mathsf{s},\mathsf{s}_{\pm}},\; \ell(\mathsf{s},y)\leq 2/3
\ell(\mathsf{s},\mathsf{s}_{\pm})\}$
where $ \Gamma_{\mathsf{s},\mathsf{s}_{\pm}}$ denotes the piece of the boundary of $\Omega$ which joins the edges $\mathsf{s}$ and $\mathsf{s}_{\pm}$ and $\ell(\mathsf{s},\mathsf{s}_{\pm})$ is the length of $\Gamma_{\mathsf{s},\mathsf{s}_{\pm}}$.
Figures~\ref{fig1} and \ref{fig2} give examples of domains $\Omega_{\mathsf{s}}$ and $\Omega_{\mathsf{s}}^\pm$. 

\begin{figure}[h!]
\begin{center}
 \includegraphics[height=4.5cm]{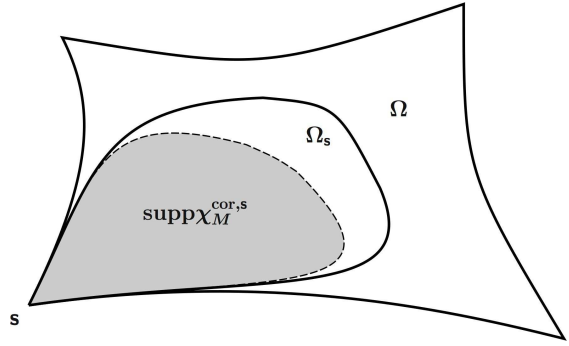}
 \caption{Definition of $\Omega_{\mathsf{s}}$}\label{fig1} 
\end{center}
\end{figure}
\begin{figure}[h!]
\begin{center}
 \includegraphics[height=4.5cm]{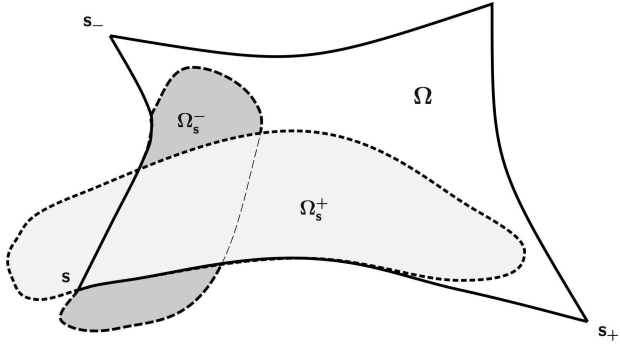}
 \caption{Definition of $\Omega_{\mathsf{s}}^+$ and $\Omega_{\mathsf{s}}^-$}\label{fig2}
\end{center}
\end{figure}

As soon as $B/M^2$ is large enough, the support of $\chi^{\rm cor}_M$ is the union of $N$ disjoint domains localized near each corner $\mathsf{s}$, $\mathsf{s}\in\Sigma$. Consequently, for $B\geq B_{0}$, we can rewrite $\chi^{\rm cor}_M$ as
$$\chi^{\rm cor}_M = \sum_{\mathsf{s}\in\Sigma} \chi^{{\rm cor},\mathsf{s}}_M\quad\mbox{ with }\;\mathsf{s}\in{\rm supp }\chi^{{\rm cor},\mathsf{s}}_M,\quad {\rm supp }\chi^{{\rm cor},\mathsf{s}}_M\cap {\rm supp }\chi^{{\rm cor},\mathsf{s}'}_M=\emptyset, \forall \mathsf{s}\neq \mathsf{s}'.$$
Furthermore, we choose $B_{0}$ large enough such that for any $B\geq B_{0}$,
$${\rm supp }\chi_M^{{\rm cor},\mathsf{s}} \cap \Omega \subset \Omega_{\mathsf{s}},\quad \forall \mathsf{s}\in\Sigma.$$
Using the eigenvalue asymptotics from \cite[Prop. 11.4]{Bo2} and \cite[Th. 7.1]{BoDa}, we therefore conclude that
\begin{align}
\label{eq:COR}
Q_{\Omega,B}[ \chi^{{\rm cor},\mathsf{s}}_M\phi]
\geq (\mu_1(\alpha_{\mathsf{s}})B - C B^{1/2}) \| \chi^{{\rm cor},\mathsf{s}}_M\phi \|^2.
\end{align}
By a similar argument, we prove an analogous lower bound for the boundary contribution. Indeed, 
if $B$ is large enough, the support of $\chi_M^{\rm bd}$ is the union of $N$ disjoint (c.f. the choice of $L$) domains localized near each piece of the smooth boundary and we rewrite
$$
[\chi^{\rm bd}_M]^2 = \sum_{\mathsf{s}\in \Sigma}([\chi^{{\rm bd},\mathsf{s},-}_M]^2+[\chi^{{\rm bd},\mathsf{s},+}_M]^2)\quad\mbox{ with }\;
 {\rm supp }\chi^{{\rm bd},\mathsf{s},\pm}_M\subset \Omega_{\mathsf{s}}^\pm\cap \Omega,\quad \forall \mathsf{s}\in\Sigma.$$
Let $\mathsf{s}\in\Sigma$. 
 From the asymptotics of the ground state energy of ${\mathcal H}_{\Omega'}(B)$ for smooth domains $\Omega'$ (\cite[Thm.~11.1]{He-Mo}) we get the following lower bound
\begin{align}
\label{eq:BDRY.0}
Q_{\Omega,B}[ \chi^{{\rm bd},\mathsf{s},\pm}_M\phi] &=
Q_{\Omega_{\mathsf{s}}^{\pm},B}[ \chi^{{\rm bd},\mathsf{s},\pm}_M\phi] \nonumber\\
&\geq 
\big\{\Theta_{0}B -2M_{3}B^{1/2}\kappa_{\rm max}^{\pm}(\mathsf{s})-C_{0}(\Omega_{\mathsf{s}}^{\pm})B^{1/3}\big\}\| \chi^{{\rm bd},\mathsf{s},\pm}_M\phi \|^2,
\end{align}
where $M_3$ is a universal constant, $C_{0}(\Omega_{\mathsf{s}}^{\pm})$ is a domain-dependent constant and $\kappa^{\pm}(\mathsf{s})$ denotes the maximal curvature of the boundary $\partial \Omega_{\mathsf{s}}^{\pm}$. We can bound $\kappa_{\rm max}(\mathsf{s})$ by 
$$
\kappa_{\rm max}:=\max_{\mathsf{s} \in\Sigma,\pm} \kappa_{\rm max}^{\pm}(\mathsf{s}),
$$
and similarly for $C_{0}(\Omega_{\mathsf{s}}^{\pm})$. 
Then, there exists $C$ independent of $M$ and $\mathsf{s}$ such that
\begin{align}
\label{eq:BDRY.1}
Q_{\Omega,B}[ \chi^{{\rm bd},\mathsf{s},\pm}_M\phi] 
\geq (\Theta_{0}B -CB^{1/2})\| \chi^{{\rm bd},\mathsf{s},\pm}_M\phi \|^2.
\end{align}
Using again the IMS-formula and \eqref{eq:BDRY.1}, we can bound
\begin{align}
\label{eq:BDRY}
Q_{\Omega,B}[ \chi^{{\rm bd}}_M\phi] &\geq
\sum_{\mathsf{s} \in \Sigma, \pm}Q_{\Omega,B}[ \chi^{{\rm bd},\mathsf{s},\pm}_M\phi]  - C \frac{B}{M^2} \| \phi \|^2
\nonumber\\
&\geq (\Theta_{0}B -CB^{1/2})\| \chi^{\rm bd}_M\phi \|^2  - C \frac{B}{M^2} \| \phi \|^2.
\end{align}
We clearly get the result of Theorem~\ref{thm:LowerBoundPotentialCorners} by combining \eqref{eq:IMS} with \eqref{eq:INT}, \eqref{eq:COR}  and \eqref{eq:BDRY} and choosing $M_0, B_0$ sufficiently large.
\end{proof}

Using the lower bound \eqref{eq:LowerInterior} combined with the upper bound inherent in \eqref{eq:EigenvalueAsymp}, one can get the following Agmon type estimate for the linear problem.
Again we only state the result for completeness and without proof, since we will not use Theorem~\ref{thm:AgmonCornersLinear} in the remainder of the paper.

\begin{thm}\label{thm:AgmonCornersLinear}~\\
Let $\psi_B$ be the ground state eigenfunction of ${\mathcal H}(B)$.
Then there exist constants $\epsilon, C, B_0>0$ such that
$$
\int e^{\epsilon \sqrt{B} \dist(x, \Sigma)}\big\{ |\psi_B(x)|^2 + B^{-1} |p_{B{\bf F}}\psi_B(x)|^2\big\}\,dx
\leq
C \| \psi_B\|_2^2\;,
$$
for all $B\geq B_0$.
\end{thm}

\section{Basic estimates}
We will need a number of standard results that we collect here for easy reference.
First of all we have the usual $L^{\infty}$-bound for solutions to the Ginzburg-Landau equations \eqref{eq:GL},
\begin{align}
\label{eq:MaxPr}
\| \psi \|_{\infty} \leq 1.
\end{align}
The proof in \cite{DGP} does not depend on regularity of the boundary, in particular, it is valid for domains with Lipschitz boundary.

The normalization of our functional ${\mathcal E}_{\kappa, H}$ is such that ${\mathcal E}_{\kappa, H}[0,{\bf F}]=0$. So any minimizer $(\psi, {\bf A})$ will have non-positive energy. Therefore, the only negative term, $-\kappa^2 \| \psi \|_2^2$, in the functional has to control each of the positive terms. This leads to the following basic inequalities for minimizers,
\begin{align}
\label{eq:BasicQuadForm}
&\| p_{\kappa H {\bf A}} \psi \|_2 \leq \kappa \|\psi \|_2, \\
\label{eq:BasicCurl}
&H \| \curl {\bf A} -1 \|_2 \leq \| \psi \|_2.
\end{align}
Furthermore, using \eqref{eq:MaxPr},
\begin{align}
\label{eq:Basic4-2}
\| \psi \|_4^2 \leq  \| \psi \|_2.
\end{align}
The following lemma states that 
in two dimensions it is actually irrelevant whether we integrate the fields over $\Omega$ or over ${\mathbb R}^2$ in the definition of ${\mathcal E}_{\kappa, H}$.

\begin{lemma}
\label{lem:EqualOutside}~\\
Let $\Omega$ be a bounded domain with Lipschitz boundary and let $(\psi, {\bf A})$ be a (weak) solution to \eqref{eq:GL}. Then $\curl({\bf A} - {\bf F}) = 0$ on the unbounded component of ${\mathbb R}^2 \setminus \overline{\Omega}$.
\end{lemma}

\begin{proof}~\\
The second equation, \eqref{eq:equationB} reads in the exterior of $\Omega$, using that $\curl {\bf F} =1$,
$$
\big(\partial_2 \curl ({\bf A} - {\bf F}), -\partial_1 \curl ({\bf A} - {\bf F})\big) = 0.
$$
Thus we see that $\curl ({\bf A} - {\bf F})$ is constant on each connected component of ${\mathbb R}^2 \setminus \overline{\Omega}$ and since it has to be in $L^2$ it must therefore vanish on the unbounded component.
\end{proof}

\begin{lemma}
\label{lem:H1ControlsL2}~\\
There exists a constant $C_0$ (depending only on $\Omega$) such that if
$(\psi, {\bf A})$ is a (weak) solution of the Ginzburg-Landau equations \eqref{eq:GL}, then
\begin{align}
\label{eq:L2again}
\int_{\Omega} | {\bf A} - {\bf F}|^2 \leq C_0 \int_{{\mathbb R}^2} | \curl {\bf A} -1 |^2 \,dx,\\
\label{eq:H1}
\| {\bf A} - {\bf F}\|^2_{W^{1,2}(\Omega)} \leq C_0 \int_{{\mathbb R}^2} | \curl {\bf A} -1 |^2 \,dx.
\end{align}
\end{lemma}

\begin{proof}~\\
Let $b=\curl({\bf A} - {\bf F})$. By Lemma~\ref{lem:EqualOutside}, $\supp b \subseteq \overline{\Omega}$. Define $\Gamma_2(x) = \frac{1}{2\pi} \log(|x|)$ (the fundamental solution of the Laplacian in two dimensions), and $w=\Gamma_2 * b$. Then $w \in \dot{H}^2({\mathbb R}^2)$ and 
(see \cite{GilbargTrudinger})
\begin{align}
\label{eq:GT}
\Delta w = b, \quad\quad
\| \nabla w \|_{L^2(\Omega)} \leq C(\Omega) \| b\|_{L^2(\Omega)}.
\end{align}
Let $\tilde{\bf A} = (-\partial_2 w, \partial_1 w) \in \dot{H}^1({\mathbb R}^2)$. Then
$$
\Div \tilde{\bf A} = 0,\quad\quad \curl \tilde{\bf A} = b.
$$
So we conclude that $\tilde{\bf A} = {\bf A} - {\bf F}$, and therefore \eqref{eq:L2again} follows from \eqref{eq:GT}.

To establish \eqref{eq:H1} we use \eqref{eq:L2again} together with the standard estimate
$$
\| D{\bf a} \|_{L^2({\mathbb R}^2)} \leq C 
\big( \|\Div {\bf a} \|_{L^2({\mathbb R}^2)}+ \|\curl {\bf a} \|_{L^2({\mathbb R}^2)}\big).
$$
\end{proof}

\section{Non-linear Agmon estimates}
\subsection{Rough bounds on $\| \psi \|_2^2$}~\\
In this chapter we prove that minimizers are localized near the boundary when $H>\kappa$. 
The precise meaning of that statement is given by Theorem~\ref{thm:Rough} below. In particular, since $\| \psi \|_{\infty} \leq 1$, the $L^2$-norm satisfies $\| \psi \|_2 = o(1)$.
We thus give a very precise and general upper bound to the field strength above which superconductivity is essentially a boundary phenomenon. Notice that this is the field which is usually called $H_{C_2}$ in the literature, although a precise mathematical definition is somewhat difficult to give.

The proof of Theorem~\ref{thm:Rough} given below has been developed in cooperation with R.~Frank.

\begin{theorem}[Weak decay estimate]
\label{thm:Rough}~\\
Let $\Omega$ be a bounded domain with Lipschitz boundary.
Then there exist $C,C'>0$, such that if  
$(\psi, {\bf A})_{\kappa,H}$ is a minimizer of ${\mathcal E}_{\kappa,H}$ with 
\begin{align}
\label{eq:AboveC2}
\kappa(H-\kappa) \geq 1/2,
\end{align}
then
\begin{align}\label{eq:RoughL2}
\| \psi \|_2^2 \leq C \int_{ \{\sqrt{\kappa(H-\kappa)}\, \dist(x,\partial \Omega)\leq 1 \}} |\psi(x)|^2\,dx 
\leq \frac{C'}{\sqrt{\kappa(H-\kappa)}}.
\end{align}
\end{theorem}

\begin{proof}~\\
The last inequality is an easy consequence of \eqref{eq:MaxPr}, since there exists a constant $C_1>0$ (depending only on $\Omega$) such that
$\meas \{ x\; : \;  \dist(x,\partial \Omega) \leq \lambda\} \leq C_1 \lambda$ for all $\lambda \in (0,2]$.

Let $\chi \in C^{\infty}({\mathbb R})$ be a standard non-decreasing cut-off function,
$$
\chi = 1 \quad \text{ on } [1,\infty), \quad\quad \chi = 0 \quad \text{ on } (-\infty,1/2).
$$
Notice for later use that this implies that $\| \chi' \|_{\infty} \geq 2$.
Let further $\lambda>0$ (we will choose $\lambda={1}/{\sqrt{\kappa(H-\kappa)}}$ at the end of the proof) and define $\chi_{\lambda}:\Omega \rightarrow {\mathbb R}$ by
$$
\chi_{\lambda}(x)  = \chi(\dist(x,\partial \Omega)/\lambda).
$$
Then $\chi_{\lambda}$ is a Lipschitz function and $\supp \chi_{\lambda} \subset \Omega$. Combining the standard localization formula and \eqref{eq:equationA}, we find
\begin{align}
\label{eq:LocFirst}
\int_{\Omega} | p_{\kappa H {\bf A}} (\chi_{\lambda} \psi) |^2 \,dx-
\int_{\Omega} |\nabla \chi_{\lambda}|^2 |\psi |^2\,dx
&=
\Re \langle \chi_{\lambda}^2 \psi,  {\mathcal H}_{\kappa H {\bf A}} \psi \rangle \nonumber\\
&=
\kappa^2 \int |\chi_{\lambda} \psi|^2\,dx - \kappa^2 \int \chi_{\lambda}^2 |\psi|^4\,dx.
\end{align}
Since $\chi_{\lambda} \psi$ has compact support we have
\begin{align}
\label{eq:CompSuppB}
\int_{\Omega} | p_{\kappa H {\bf A}} (\chi_{\lambda} \psi) |^2\,dx
&\geq
\kappa H \int_{\Omega} (\curl {\bf A}) |\chi_{\lambda}\psi |^2 \,dx\nonumber\\
&\geq
\kappa H \|  \chi_{\lambda} \psi \|_2^2  - 
\kappa H \| \curl {\bf A} -1 \|_2 \|  \chi_{\lambda} \psi \|_4^2.
\end{align}
Using \eqref{eq:Basic4-2} and \eqref{eq:BasicCurl}, we get from \eqref{eq:LocFirst} and \eqref{eq:CompSuppB} that
\begin{align*}
&\kappa (H- \kappa) \|  \chi_{\lambda} \psi \|_2^2\\
&\leq 
\kappa \| \psi \|_2 \| \chi_{\lambda} \psi \|_4^2
- \kappa^2 \int \chi_{\lambda}^2 |\psi|^4\,dx + \| \chi' \|_{\infty}^2 \lambda^{-2} \int_{\{\dist(x,\partial \Omega) \leq \lambda\}}
|\psi(x)|^2\,dx \\
&\leq
\frac{1}{4} \| \psi \|_2^2 + \| \chi' \|_{\infty}^2 \lambda^{-2} \int_{\{\dist(x,\partial \Omega) \leq \lambda\}}
|\psi(x)|^2\,dx  + \kappa^2 \int (\chi_{\lambda}^4- \chi_{\lambda}^2) |\psi|^4\,dx.
\end{align*}
Notice that since $\chi\leq 1$, the last integral is negative and we thus find by 
dividing the integral $\| \psi \|_2^2$ in two
\begin{align*}
\{\kappa (H- \kappa) -&1/4 \} \|  \chi_{\lambda} \psi \|_2^2\\
&\leq
\frac{1}{4} \int (1 - \chi_{\lambda}^2) | \psi |^2\,dx
+ \| \chi' \|_{\infty}^2 \lambda^{-2} \int_{\{\dist(x,\partial \Omega) \leq \lambda\}}
|\psi(x)|^2\,dx \\
&\leq (\| \chi' \|_{\infty}^2 \lambda^{-2} +1/4) \int_{\{\dist(x,\partial \Omega) \leq \lambda\}}
|\psi(x)|^2\,dx.
\end{align*}
Choose $\lambda = 1/{ \sqrt{\kappa(H-\kappa)}}$. By assumption 
$\kappa(H-\kappa) -1/4 \geq \kappa(H-\kappa)/2$, and the conditions on $\chi$, $\kappa(H-\kappa)$ imply that
$ \| \chi' \|_{\infty}^2 \lambda^{-2} +1/4 \leq 2 \| \chi' \|_{\infty}^2 \lambda^{-2}$.
Thus,
\begin{align}
\|  \chi_{\lambda} \psi \|_2^2
\leq 4 \| \chi' \|_{\infty}^2 \int_{\{\dist(x,\partial \Omega) \leq \lambda\}}
|\psi(x)|^2\,dx.
\end{align}
Consequently,
\begin{align}
\|  \psi \|_2^2
\leq (4 \| \chi' \|_{\infty}^2 + 1)\int_{\{\dist(x,\partial \Omega) \leq \lambda\}}
|\psi(x)|^2\,dx.
\end{align}
This finishes the proof of \eqref{eq:RoughL2}.
\end{proof}

For stronger fields superconductivity is essentially localized to the corners.

\begin{theorem}[Decay estimate on the boundary]
\label{thm:Rough2}~\\
Suppose that $\Omega$ satisfies Assumption~\ref{ass:Angles}.
For $\mu \in ( \Lambda_{1}, \Theta_0)$, define
\begin{align}
\Sigma':=\{\mathsf{s}\in\Sigma\, \big |\, \mu_{1}(\alpha_{\mathsf{s}})\leq \mu\}, \quad\text{ and }\quad
b:=\inf_{\mathsf{s} \in \Sigma \setminus \Sigma'} \{ \mu_{1}(\alpha_{\mathsf{s}}) - \mu\}.
\end{align}
(in the case $\Sigma = \Sigma'$, we set $b:=\Theta_0 - \mu$).\\
There exist $\kappa_0,C,C',M>0$, such that if $(\psi, {\bf A})_{\kappa,H}$ is a minimizer of ${\mathcal E}_{\kappa,H}$ with 
\begin{align}
\label{eq:AboveC2.2}
\frac{H}{\kappa} \geq \mu^ {-1}, \quad\quad \kappa \geq \kappa_0,
\end{align}
then
\begin{align}\label{eq:RoughL2.2}
\| \psi \|_2^2 \leq C \int_{ \{\kappa \, \dist(x,{\Sigma'  })\leq M \}} |\psi(x)|^2\,dx 
\leq \frac{C'}{\kappa^2}.
\end{align}
\end{theorem}

\begin{proof}~\\
To prove this result, we follow the same procedure as in the proof of Theorem~\ref{thm:Rough}.\\
Let $\delta=b/2$, and let $M_0=M_0(\delta)$ be the constant from Theorem~\ref{thm:LowerBoundPotentialCorners}.
Let $\chi \in C^{\infty}({\mathbb R})$ be a standard non-decreasing cut-off function,
$$
\chi = 1 \quad \text{ on } [1,\infty), \quad\quad \chi = 0 \quad \text{ on } (-\infty,1/2),
$$
and let $\lambda={2M_0}/{\sqrt{\kappa H}}$.
Define $\chi_{\lambda}:\Omega \rightarrow {\mathbb R}$, by
$$
\chi_{\lambda}(x)  = \chi(\dist(x,{\Sigma'})/\lambda).
$$
Then $\chi_{\lambda}$ is a Lipschitz function and $\supp \chi_{\lambda} \cap {\Sigma'}=\emptyset$. Combining the standard localization formula and \eqref{eq:equationA}, we find as previously
\begin{align}
\label{eq:LocFirst.2}
\int_{\Omega} | p_{\kappa H {\bf A}} (\chi_{\lambda} \psi) |^2 \,dx-
\int_{\Omega} |\nabla \chi_{\lambda}|^2 |\psi |^2\,dx
=\Re \langle \chi_{\lambda}^2 \psi,  {\mathcal H}_{\kappa H {\bf A}} \psi \rangle 
&\leq \kappa^2 \| \chi_{\lambda} \psi \|_2^2.
\end{align}
As in \eqref{eq:CompSuppB}, we need a lower bound to $\int_{\Omega} | p_{\kappa H {\bf A}} (\chi_{\lambda} \psi) |^2\,dx$. Since $\supp \chi_{\lambda} \cap \partial\Omega\neq\emptyset$, 
we cannot argue as in \eqref{eq:CompSuppB}. Therefore,
we will introduce the constant magnetic field ${\bf F}$ for which we have such an estimate, namely Theorem~\ref{thm:LowerBoundPotentialCorners}. We can write
\begin{align}
\label{eq:CompSuppB.2}
\int_{\Omega} | p_{\kappa H {\bf A}} (\chi_{\lambda} \psi) |^2\,dx
&\geq (1-\varepsilon)\int_{\Omega} | p_{\kappa H {\bf F}} (\chi_{\lambda} \psi) |^2\,dx \nonumber\\
&\quad\quad
-\varepsilon^{-1} \int_{\Omega} (\kappa H)^2 |{\bf F}-{\bf A} |^2  (\chi_{\lambda} \psi) |^2\,dx. 
\end{align}
Theorem~\ref{thm:LowerBoundPotentialCorners} and the choice of $\lambda$ imply that
\begin{align}
\int_{\Omega} | p_{\kappa H {\bf F}} (\chi_{\lambda} \psi) |^2\,dx 
&\geq \left(\inf_{\mathsf{s}\in\Sigma\setminus\Sigma'}\mu_{1}(\alpha_{\mathsf{s}})-\delta\right)
\kappa H \|  \chi_{\lambda} \psi \|_2^2 \nonumber\\
&= \left({\mu}+\frac b2\right) \kappa H \|  \chi_{\lambda} \psi \|_2^2.\label{eq:CompSuppB.3}
\end{align}
We now have to give a lower bound to the second part of the right side of \eqref{eq:CompSuppB.2}. We can estimate 
\begin{equation}\label{eq:CompSuppB.4}
\int_{\Omega} (\kappa H)^2 |{\bf F}-{\bf A} |^2  |\chi_{\lambda} \psi |^2\,dx
\leq (\kappa H)^2 \|{\bf A}-{\bf F}\|^2_{4}\; \|\chi_{\lambda}\psi\|^2_{4}.
\end{equation}
By Sobolev inequalities, \eqref{eq:H1} and \eqref{eq:BasicCurl}, we deduce
\begin{align}\label{eq:CompSuppB.5}
(\kappa H)^2 \|{\bf F}-{\bf A} \|^2_{4} 
&\leq C \kappa^2 H^2\|{\bf F}-{\bf A} \|^2_{W^{1,2}(\Omega)} \nonumber\\
&\leq  \tilde C \kappa^2 H^2 \|{\rm curl}{\bf A}-1\|^2_{L^2(\R^2)}\nonumber\\
&\leq \tilde C \kappa^2 \|\psi\|^2_{2}.
\end{align}
Let us now estimate $\|\chi_{\lambda}\psi\|^2_{4}$.  According to \eqref{eq:MaxPr} and the property of the cut-off function $0\leq\chi_{\lambda}\leq 1$, we can bound $|\chi_{\lambda}\psi|$ from above by $1$ and deduce, using also Theorem~\ref{thm:Rough},
\begin{equation}\label{eq:CompSuppB.6}
\|\chi_{\lambda}\psi\|^2_{4} = \sqrt{\int_{\Omega} |\chi_{\lambda}\psi|^4\, dx}
\leq \sqrt{\int_{\Omega} |\chi_{\lambda}\psi|^2\, dx}
\leq \frac C{\sqrt \kappa}.
\end{equation}
Inserting \eqref{eq:CompSuppB.3}, \eqref{eq:CompSuppB.4}, \eqref{eq:CompSuppB.5} and \eqref{eq:CompSuppB.6} in \eqref{eq:CompSuppB.2}, we obtain
\begin{equation}\label{eq:CompSuppB.7}
\int_{\Omega} | p_{\kappa H {\bf A}} (\chi_{\lambda} \psi) |^2\,dx
\geq  (1-\varepsilon)\left({\mu} + \frac b2 \right)\kappa H  \|  \chi_{\lambda} \psi \|_2^2 
-C\varepsilon^{-1}  \kappa^{3/2} \|\psi\|^2_{2}.
\end{equation}
We insert \eqref{eq:CompSuppB.7} in \eqref{eq:LocFirst.2}. Then
\begin{multline}\label{eq:CompSuppB.80}
\left[(1-\varepsilon)\left(\mu+\frac b2\right)\kappa H-\kappa^2-C\varepsilon^{-1}\kappa^{3/2}\right]
\int_{\{{\rm dist}(x,{\Sigma'})\geq \lambda\}}| \psi|^2\, dx\\
\leq 
(C\varepsilon^{-1}\kappa^{3/2}+\|\chi'\|_{\infty}^2\lambda^{-2})\int_{\{{\rm dist}(x,{\Sigma'})\leq\lambda\}}|\psi|^2\, dx,
\end{multline}
Assumption~\eqref{eq:AboveC2.2} leads to the lower bound
\begin{align}\label{eq:CompSuppB.9}
(1-\varepsilon)\left(\mu+\frac b2\right)\kappa H-\kappa^2-C\varepsilon^{-1}\kappa^{3/2}
\geq
\frac{b}{4} \kappa H,
\end{align}
as soon $\varepsilon$ is small enough and $\kappa$ large enough.

Once $\varepsilon$ is fixed and with $\lambda = {2M_0}/{\sqrt{\kappa H}}$, we find
\begin{equation}\label{eq:CompSuppB.10}
C\varepsilon^{-1}  \kappa^{3/2} +\|\chi'\|_{\infty}\lambda^{-2} \leq c \kappa H.
\end{equation}
Combining \eqref{eq:CompSuppB.80}, \eqref{eq:CompSuppB.9} and \eqref{eq:CompSuppB.10}, we deduce
\begin{equation}\label{eq:CompSuppB.11}
\int_{\{{\rm dist}(x,{\Sigma'})\geq \lambda\}}| \psi|^2\, dx
\leq C \int_{\{{\rm dist}(x,{\Sigma'})\leq \lambda\}}| \psi|^2\, dx.
\end{equation}
It follows easily that 
$$\| \psi\|_{2}^2 \leq (C+1)  \int_{\{{\rm dist}(x,{\Sigma'})\leq \lambda\}}| \psi|^2\, dx.
$$
Inserting the choice $\lambda = {2M_0}/{\sqrt{\kappa H}}$ and the condition \eqref{eq:AboveC2.2} on $H$, this clearly implies \eqref{eq:RoughL2.2}.
\end{proof}

\subsection{Exponential localization}~\\
In order to obtain exponential decay in the interior of the domain, we need the following energy estimate, Lemma~\ref{lem:EnergyInterior}, for functions located away from the boundary.

\begin{lemma}\label{lem:EnergyInterior}~\\
Let $\Omega \subset {\mathbb R}^2$ be a bounded domain with Lipschitz boundary.
There exist constants $C_0, C_1>0$ such that if $\kappa(H-\kappa) \geq C_0$ and $(\psi, {\bf A})$ 
is a minimizer of ${\mathcal E}_{\kappa,H}$, then for all $\phi \in C_0^{\infty}(\Omega)$ we have
$$
\| (-i\nabla - \kappa H {\bf A}) \phi \|_{2}^2 \geq \kappa H \big (1 - C_1 \| \psi \|_2 \big )
\| \phi \|_{2}^2.
$$
In particular, using the estimate on $\| \psi \|_2$ from Theorem~\ref{thm:Rough} we find
$$
\| (-i\nabla - \kappa H {\bf A}) \phi \|_{2}^2 \geq \kappa H \left(1 - \frac{C_1'}{\sqrt[4]{\kappa(H-\kappa)}} \right) \| \phi \|_{2}^2.
$$
\end{lemma}

\begin{proof}~\\
We estimate, for $\phi \in C_0^{\infty}(\Omega)$, 
\begin{align}
\label{eq:EnergyComp}
\| (-i\nabla - \kappa H {\bf A}) \phi \|_{2}^2 &\geq
\kappa H \int_{\Omega} \curl {\bf A} |\phi|^2\,dx \nonumber \\
&\geq \kappa H \| \phi \|_2^2 - \kappa H \| \curl {\bf A} - 1 \|_2 \| \phi \|_4^2.
\end{align}
By the Sobolev inequality, for $\phi \in C_0^{\infty}({\mathbb R}^2)$, and scaling we get, for all $\eta >0$ and with a universal constant $C_{\rm Sob}$, the estimate
\begin{align}
\label{eq.Phi4}
\| \phi \|_4^2 \leq C_{\rm Sob}\Big( \eta \big \| \nabla|\phi| \big\|_2^2 + \eta^{-1} \| \phi \|_2^2 \Big).
\end{align}
We can estimate $\| \nabla|\phi| \big\|_2^2$ by $\| (-i\nabla - \kappa H {\bf A}) \phi \|_2^2$ by the diamagnetic inequality. Choosing, $\eta = \frac{\eta'}{C_{\rm Sob} \kappa H \| \curl {\bf A} - 1 \|_2}$, for some $\eta'>0$, we thus find, using \eqref{eq:BasicCurl}, \eqref{eq:EnergyComp} and \eqref{eq.Phi4},
\begin{align}
\| (-i\nabla& - \kappa H {\bf A}) \phi \|_2^2 \nonumber\\
&\geq \kappa H \| \phi \|_2^2 - \eta' \| (-i\nabla - \kappa H {\bf A}) \phi \|_2^2
- (\eta')^{-1} C^2_{\rm Sob} (\kappa H)^2 \| \curl {\bf A} - 1 \|_2^2  \| \phi \|_2^2\nonumber\\
&\geq \kappa H \| \phi \|_2^2 \Big( 1 - (\eta')^{ -1} C^2_{\rm Sob} \frac{\kappa}{H}  \| \psi \|_2^2 \Big) - \eta' \| (-i\nabla - \kappa H {\bf A}) \phi \|_2^2.
\end{align}
By assumption $\kappa/H\leq 1$. We take $\eta' = \| \psi \|_2$ and find
\begin{align}
(1+ \| \psi \|_2) \| (-i\nabla& - \kappa H {\bf A}) \phi \|_2^2
\geq \kappa H (1- C^2_{\rm Sob} \| \psi \|_2) \| \phi \|_2^2.
\end{align}
By Theorem~\ref{thm:Rough} we have
$$
\frac{1- C^2_{\rm Sob} \| \psi \|_2}{1+ \| \psi \|_2} \geq 1 - 2C^2_{\rm Sob} \| \psi \|_2,
$$
if $\kappa(H-\kappa)$ is sufficiently big. This finishes the proof of Lemma~\ref{lem:EnergyInterior}.
\end{proof}

By standard arguments Lemma~\ref{lem:EnergyInterior} implies Agmon estimates in the interior.

\begin{thm}[Normal Agmon estimates]\label{thm:AgmonNormal}~\\
Let $\Omega$ be a bounded domain with Lipschitz boundary and let $b>0$. There exist $M,C,\epsilon, \kappa_0>0$, such that if $(\psi, {\bf A})$ is a minimizer of ${\mathcal E}_{\kappa,H}$ with 
$$
\frac{H}{\kappa} \geq 1 + b,\quad \kappa \geq \kappa_0,
$$ 
then
\begin{align}
\label{eq:AgmonNormal}
\int_{\Omega} e^{2\epsilon \sqrt{\kappa H} t(x)} \Big( |\psi|^2 + \frac{1}{\kappa H} \big| (-i\nabla - \kappa H {\bf A})\psi \big|^2 \Big)\,dx \leq C \int_{\{ t(x)\leq \frac{M}{\sqrt{\kappa H}}\}} |\psi |^2\,dx.
\end{align}
Here $t(x) := \dist(x, \partial \Omega)$.
\end{thm}

\begin{proof}~\\
The function $t(x) = \dist(x, \partial \Omega)$ defines a Lipschitz continuous function on $\Omega$. In particular, $\nabla t \in L^{\infty}(\Omega)$. Let $\chi \in C^{\infty}({\mathbb R})$ be a non-decreasing function satisfying
$$
\chi = 1 \quad \text{ on } [1,\infty), \quad\quad
\chi = 0 \quad \text{ on } [-\infty, 1/2).
$$
Define the (Lipschitz continuous) function $\chi_M$ on $\Omega$ by
$\chi_M(x) = \chi(\frac{t(x) \sqrt{\kappa H}}{M})$. We calculate, using \eqref{eq:equationA} and the IMS-formula
\begin{align}
\label{eq:IMS2}
\kappa^2 \| \exp\big(\epsilon &\sqrt{\kappa H} t\big) \chi_M \psi \|_2^2
\geq
\Re \big\langle \exp\big(2\epsilon \sqrt{\kappa H} t\big) \chi_M^2 \psi, \kappa^2(1-|\psi|^2)\psi \big\rangle \nonumber \\
&= \int_{\Omega} \big| p_{\kappa H {\bf A}} \big( e^{\epsilon \sqrt{\kappa H} t} \chi_M \psi\big) \big|^2\,dx - \int_{\Omega} \big| \nabla (e^{\epsilon \sqrt{\kappa H} t} \chi_M ) \psi \big|^2\,dx .
\end{align}
Combining Theorem~\ref{thm:Rough} with Lemma~\ref{lem:EnergyInterior} there exists $\tilde{g}$ with $\tilde{g} = o(1)$ at $\infty$, such that
$$
\int_{\Omega} \big| p_{\kappa H {\bf A}} \big( e^{\epsilon \sqrt{\kappa H} t} \chi_M \psi\big) \big|^2\,dx
\geq \kappa H (1 + \tilde{g}(\kappa H)) \| e^{\epsilon \sqrt{\kappa H} t} \chi_M \psi \|_2^2.
$$
Since $\frac{H}{\kappa} \geq 1 + b$, we therefore find, with some constant $C$ independent of $\kappa, H, \epsilon$ and $M$
\begin{align}
\left( 1 + \tilde{g}(\kappa H) - \frac{1}{1+b}\right) 
&\| e^{\epsilon \sqrt{\kappa H} t} \chi_M \psi \|_2^2 \nonumber\\
&\leq
C \epsilon^2 \| \nabla t \|_{\infty}^2 \| e^{\epsilon \sqrt{\kappa H} t} \chi_M \psi \|_2^2\nonumber\\
&\quad
+
\frac{C \| \nabla t \|_{\infty}^2}{M^2} \int_{\Omega} e^{2 \epsilon \sqrt{\kappa H} t(x)} \bigg| \chi'\bigg(\frac{t(x)\sqrt{\kappa H}}{M}\bigg) \psi(x)\bigg|^2\,dx.
\end{align}
For $\kappa$ sufficiently big we have, since $H \geq (1+b)\kappa$, 
$$
1 + \tilde{g}(\kappa H) - \frac{1}{1+b} \geq b/2.
$$
We choose $\epsilon$ sufficiently small that $C \epsilon^2 \| \nabla t \|_{\infty}^2 < b/4$ and finally obtain for some new constant $C'$
\begin{align}
\| e^{\epsilon \sqrt{\kappa H} t} \chi_M \psi \|_2^2 
\leq C' \frac{e^{2\epsilon M}}{M^2} \int_{\{\sqrt{\kappa H} t(x)\leq M\}} |\psi(x)|^2\,dx.
\end{align}
On the support of $1-\chi_M$ the exponential $e^{\epsilon \sqrt{\kappa H} t}$ is bounded, so we see that 
\begin{align}
\label{AgmonL2}
\| e^{\epsilon \sqrt{\kappa H} t}  \psi \|_2^2 
\leq C'' \int_{\{\sqrt{\kappa H} t(x)\leq M\}} |\psi(x)|^2\,dx,
\end{align}
which is part of the estimate \eqref{eq:AgmonNormal}. 

It remains to estimate the term with $\big| (-i\nabla - \kappa H {\bf A})\psi \big|$ in \eqref{eq:AgmonNormal}. This follows from the same considerations upon inserting the bound \eqref{AgmonL2}.
\end{proof}

\begin{lemma}\label{lem:EnergyBordNew}~\\
Suppose that $\Omega \subset {\mathbb R}^2$ satisfies Assumption~\ref{ass:Angles}.
For $\mu \in ( \Lambda_{1}, \Theta_0)$, define
\begin{align}
\Sigma':=\{\mathsf{s}\in\Sigma\,\big |\, \mu_{1}(\alpha_{\mathsf{s}})\leq \mu\}, \quad\text{ and }\quad
b:=\inf_{\mathsf{s} \in \Sigma \setminus \Sigma'} \{ \mu_{1}(\alpha_{\mathsf{s}}) - \mu\}.
\end{align}
(in the case $\Sigma = \Sigma'$, we set $b:=\Theta_0 - \mu$).\\
There exist $M_{0}>0$ such that if $(\psi, {\bf A})$ is a minimizer of ${\mathcal E}_{\kappa,H}$, then for all $\phi \in C^{\infty}(\overline\Omega)$ such that ${\rm dist}(\supp \phi,\Sigma')\geq M_{0}/\sqrt{\kappa H}$, we have
\begin{align}
\label{eq:LowerFormCorner2}
\| (-i\nabla-\kappa H {\bf A}) \phi \|_{L^2(\Omega)}^2 
\geq \mu\kappa H \left(1+\frac b4\right) \| \phi \|_{L^2(\Omega)}^2,
\end{align}
for $\kappa H$ sufficiently large.
\end{lemma}

\begin{proof}~\\
Let $\delta=b/2$ and let $M_0 = M_0(\delta)$ be the constant from Theorem~\ref{thm:LowerBoundPotentialCorners}.
We estimate, for $\phi \in C^{\infty}(\overline\Omega)$ such that ${\rm dist}(\supp \phi,\Sigma')\geq M_{0}/\sqrt{\kappa H}$, 
\begin{multline}
\label{eq:EnergyComp.2-2}
\| (-i\nabla - \kappa H {\bf A}) \phi \|_{2}^2 
=  \| (-i\nabla - \kappa H {\bf F}) \phi +\kappa H({\bf F}-{\bf A})\phi \|_{2}^2  \\
\geq (1-\varepsilon)\int_{\Omega}|(-i\nabla-\kappa H{\bf F})\phi|^2\, dx 
-\varepsilon^{-1} \int_{\Omega}(\kappa H)^2 |{\bf F}-{\bf A}|^2\, |\phi|^2\, dx.
\end{multline}
Using Theorem~\ref{thm:LowerBoundPotentialCorners}
and the support properties of $\phi$, we have
\begin{align}\label{eq:EnergyComp.3-2}
\int_{\Omega}|(-i\nabla-\kappa H{\bf F})\phi|^2\, dx
&\geq \left(\inf_{\mathsf{s}\in\Sigma\setminus\Sigma'}\mu_{1}(\alpha_{\mathsf{s}})-\delta\right)
\kappa H \|\phi\|^2_{2} \nonumber\\
&=\left(\mu+\frac b2\right) \kappa H \|\phi\|^2_{2}.
\end{align}
Using the Cauchy-Schwarz inequality, \eqref{eq:CompSuppB.5} and Theorem~\ref{thm:Rough2}, we can bound the last term of \eqref{eq:EnergyComp.2-2}.
\begin{align}
\label{eq:EnergyComp.4-2}
\int_{\Omega}(\kappa H)^2 |{\bf F}-{\bf A}|^2\, |\phi|^2\, dx
&\leq (\kappa H)^2 \|{\bf A}-{\bf F}\|^2_{4}\, \|\phi\|^2_{4}\nonumber\\
&\leq C \kappa^2 \|\psi\|^2_{2} \|\phi\|^2_{4}\nonumber\\
&\leq \tilde C \big\| |\phi| \big\|^2_{4}.
\end{align}
We use the Sobolev inequality \eqref{eq.Phi4} in \eqref{eq:EnergyComp.4-2} and
estimate $\| \nabla|\phi| \big\|_2^2$, using the diamagnetic inequality, by $\| (-i\nabla - \kappa H {\bf A}) \phi \|_2^2$  
to obtain
\begin{align}
\label{eq:EnergyComp.5-2}
\int_{\Omega}(\kappa H)^2 |{\bf F}-{\bf A}|^2\, |\phi|^2\, dx
&\leq  C_{\rm Sob} \big(\eta \|(-i\nabla -\kappa H{\bf A})\phi\|^2_{2}+\eta^{-1}\|\phi\|^2_{2}\big).
\end{align}
Combining \eqref{eq:EnergyComp.2-2}, \eqref{eq:EnergyComp.3-2} and \eqref{eq:EnergyComp.5-2}, we deduce that 
\begin{equation}
\label{eq:EnergyComp.6-2}
\left(1+\frac{C_{\rm Sob}\eta}{\varepsilon}\right) \| (-i\nabla - \kappa H {\bf A}) \phi \|_{2}^2 
\geq \left\{(1-\varepsilon)\left(\mu+\frac b2\right) \kappa H-\frac{C_{\rm Sob}}{\varepsilon \eta}\right\} \|\phi\|^2_{2}.
\end{equation}
We choose $\eta=\frac{C_{\rm Sob}}{\varepsilon^2\kappa H}$, then \eqref{eq:EnergyComp.6-2} becomes
\begin{equation}
\label{eq:EnergyComp.7-2}
\left(1+\frac{C_{\rm Sob}^2}{\varepsilon^3 \kappa H}\right)\|(-i\nabla-\kappa H{\bf A})\phi\|_{2}^2 \geq  \kappa H \left\{(1-\varepsilon)\left(\mu+\frac b2\right) -\varepsilon \right\} \|\phi\|^2_{2}.
\end{equation}
If we choose $\varepsilon$ sufficiently small and independent of $\kappa, H$ 
(actually, since $\mu + b/2\leq 1$, $\varepsilon = b/8$ will do) 
then
\eqref{eq:LowerFormCorner2} follows.
\end{proof}

By standard arguments Lemma~\ref{lem:EnergyBordNew} implies the Agmon estimates given in Theorem~\ref{thm:AgmonCorners}.

\begin{proof}[Proof of Theorem~\ref{thm:AgmonCorners}]~\\
The function $t'(x) := \dist(x, \Sigma')$ defines a Lipschitz continuous function on $\Omega$. In particular, $|\nabla t' |\leq 1$. Let $\chi \in C^{\infty}({\mathbb R})$ be a non-decreasing function satisfying
$$
\chi = 1 \quad \text{ on } [1,\infty), \quad\quad
\chi = 0 \quad \text{ on } [-\infty, 1/2).
$$
Define the function $\chi_M$ on $\Omega$ by
$\chi_M(x) = \chi(\frac{t'(x) \sqrt{\kappa H}}{M})$.
Using Lemma~\ref{lem:EnergyBordNew} there exists $\beta>0$, such that if $M, \kappa H$ are sufficiently large, then 
$$
\int_{\Omega} \big| p_{\kappa H {\bf A}} \big( e^{\epsilon \sqrt{\kappa H} t'} \chi_M \psi\big) \big|^2\,dx
\geq \mu \kappa H (1 + \beta) \| e^{\epsilon \sqrt{\kappa H} t'} \chi_M \psi \|_2^2.
$$

Using \eqref{eq:IMS2} and the assumption $\frac{H}{\kappa} \geq \mu^{-1}$, there exists some constant $C$ independent of $\kappa, H, \epsilon$ and $M$ such that
\begin{align}
\beta \mu \| e^{\epsilon \sqrt{\kappa H} t'} \chi_M \psi \|_2^2 
&\leq
C \epsilon^2 \| \nabla t' \|_{\infty}^2 \| e^{\epsilon \sqrt{\kappa H} t'} \chi_M \psi \|_2^2\\
&\quad+
\frac{C \| \nabla t' \|_{\infty}^2}{M^2} \int_{\Omega} e^{2 \epsilon \sqrt{\kappa H} t'(x)} \bigg| \chi'\bigg(\frac{t'(x)\sqrt{\kappa H}}{M}\bigg) \psi(x)\bigg|^2\,dx.\nonumber
\end{align}
We achieve the proof of Theorem~\ref{thm:AgmonCorners} with arguments similar to the ones of the proof of Theorem~\ref{thm:AgmonNormal}.
\end{proof}

\section{Proof of Theorem~\ref{thm:EqualHC3}}
Combining Proposition~\ref{prop.lambda1} and Lemma~\ref{lem:Asymptotics} it only remains to prove \eqref{eq:EqualFields}.
We will prove that for large $\kappa$ the following two statements are equivalent.
\begin{enumerate}
\item \label{(1)} There exists a minimizer $(\psi,{\bf A})$ of ${\mathcal E}_{\kappa,H}$ with $\| \psi\|_2 \neq 0$.
\item \label{(2)} The parameters $\kappa, H$ satisfy
\begin{align}
\label{eq:HtoSmall}
\kappa^2 - \lambda_1(\kappa H) > 0.
\end{align}
\end{enumerate}
Suppose first that \eqref{eq:HtoSmall} is satisfied. Let $u_1(\kappa H)$ be the normalized ground state eigenfunction of ${\mathcal H}(\kappa H)$ and let $t>0$. Then, for $t^2 < 2 \frac{\kappa^2 - \lambda_1(\kappa H)}{\kappa^2 \| u_1(\kappa H) \|_{4}^4}$,
\begin{align}
{\mathcal E}_{\kappa, H}[t u_1(\kappa H), {\bf F}] = t^2 [\lambda_1(\kappa H) - \kappa^2] +
\frac{\kappa^2}{2} t^4 \| u_1(\kappa H) \|_{4}^4 < 0. 
\end{align}
This shows that (\ref{(2)}) implies (\ref{(1)}). 

Notice that this first part did not need the assumption that $\kappa$ is large. However, for large $\kappa$ we know that \eqref{eq:HtoSmall} is satisfied iff $H < H^{\rm lin}_{C_3}(\kappa)$ (defined in Lemma~\ref{lem:Asymptotics}).

Suppose that $(\psi,{\bf A})$ is a non-trivial minimizer of ${\mathcal E}_{\kappa,H}$. We may assume that $H > (1+b)\kappa$ for some $b>0$, because by Proposition~\ref{prop.lambda1}, \eqref{eq:HtoSmall}  is satisfied for $\kappa \geq \kappa_0$, $H < H^{\rm lin}_{C_3}(\kappa)$, where $H^{\rm lin}_{C_3}(\kappa)$ has the asymptotics given in Lemma~\ref{lem:Asymptotics}.
Furthermore, we may assume that $H \leq T \kappa$ for some $T>0$. This follows from \cite{Giorgi-Phillips}---we give the details for completeness: \\
Since $\psi\neq 0$, we have
\begin{align*}
0 &< \lambda_1(\kappa H) \| \psi \|_2^2 \leq \int_{\Omega} |p_{\kappa H {\bf F}} \psi |^2 \,dx\\
&\leq 2 \int_{\Omega}  |p_{\kappa H {\bf A}} \psi |^2 \,dx + 2 (\kappa H)^2 \int_{\Omega}  |{\bf A} - {\bf F}|^2 |\psi |^2 \,dx.
\end{align*}
We now use, \eqref{eq:MaxPr} and Lemma~\ref{lem:H1ControlsL2} to obtain
\begin{align*}
0 < \lambda_1(\kappa H) \| \psi \|_2^2 &\leq C \Big\{ \int_{\Omega}  |p_{\kappa H {\bf A}} \psi |^2 \,dx + 
(\kappa H)^2 \int_{{\mathbb R}^2} |\curl {\bf A} -1 |^2\,dx\Big\}\\
&\leq C \kappa^2 \| \psi \|_2^2,
\end{align*}
where the last inequality holds since ${\mathcal E}_{\kappa,H}[\psi, {\bf A}] \leq 0$. Since $\lambda_1(B)$ increases linearly in $B$ we deduce that $H = {\mathcal O}(\kappa)$.

From the discussion above, we know that we may assume
$$
(1+b)\kappa \leq H \leq b^{-1} \kappa,
$$
for some $b>0$.
By Theorem~\ref{thm:Rough} we therefore find, for some $C>0$,
\begin{align}
\label{eq:L2-L4}
\| \psi \|_2^2 \leq C \Big\{ \int_{\{\dist(x,\partial \Omega) \leq \frac{1}{\kappa}\}} dx\Big\}^{1/2} \| \psi \|_4^2 \leq  C' \frac{\| \psi \|_4^2}{\sqrt{\kappa}}.
\end{align}
Since $(\psi, {\bf A})$ is a non-trivial minimizer, ${\mathcal E}_{\kappa, H}[\psi, {\bf A}] \leq 0$. So we also have
\begin{align}
\label{eq:Delta}
0< \frac{\kappa^2}{2} \| \psi \|_4^4 \leq \kappa^2 \| \psi \|_2^2 - \int_{\Omega} \big| (-i\nabla - \kappa H {\bf A}) \psi \big|^2\,dx =: \Delta.
\end{align}
The inequality \eqref{eq:L2-L4} therefore becomes,
\begin{align}
\label{eq:L2-Delta}
\| \psi \|_2^2 \leq C'' \sqrt{\Delta} \kappa^{-3/2}.
\end{align}
By Cauchy-Schwarz we can estimate
\begin{align}
0 < \Delta &= \kappa^2 \| \psi \|_2^2  - \int_{\Omega} \big| \big((-i\nabla - \kappa H {\bf F}) + \kappa H ({\bf F}- {\bf A})\big)\psi \big|^2\,dx\nonumber\\
&\leq \kappa^2 \| \psi \|_2^2 - (1 - \sqrt{\Delta} \kappa^{-3/4}) \lambda_1(\kappa H) \| \psi \|_2^2\nonumber\\
&\quad+ \frac{1}{\sqrt{\Delta} \kappa^{-3/4}}(\kappa H)^2 \int_{\Omega}  |{\bf F} - {\bf A}|^2 |\psi |^2\,dx.
\end{align}
So we find, by inserting \eqref{eq:L2-Delta}, \eqref{eq:Delta} and using Cauchy-Schwarz,
\begin{align}
\label{eq:BefBootstrap}
0 < \Delta &\leq  \big(\kappa^2 -\lambda_1(\kappa H)\big) \| \psi \|_2^2
+ C'' \frac{\lambda_1(\kappa H) \sqrt{\Delta}}{\kappa^{3/4}}\kappa^{-3/2}  \sqrt{\Delta}\nonumber\\
&\quad+\frac{\kappa^{3/4}}{\sqrt{\Delta}} (\kappa H)^2 \| {\bf F} - {\bf A} \|_{4}^2 \sqrt{\frac{2\Delta}{\kappa^2}}.
\end{align}
Since ${\mathcal E}_{\kappa, H}[\psi, {\bf A}] \leq 0$, we get using Lemma~\ref{lem:H1ControlsL2} and a Sobolev imbedding,
$$
(\kappa H)^2 \| {\bf F} - {\bf A} \|_{4}^2 \leq
C (\kappa H)^2 \| \curl {\bf A} -1 \|_{L^2({\mathbb R}^2)}^2 \leq C \Delta.
$$
Inserting this in \eqref{eq:BefBootstrap} yields,
$$
0 < \Delta \leq \big(\kappa^2 -\lambda_1(\kappa H)\big) \| \psi \|_2^2 + C \frac{\Delta}{\kappa^{1/4}}, 
$$
which permits to conclude that \eqref{eq:HtoSmall} is satisfied.

Thus (\ref{(1)}) and (\ref{(2)}) are equivalent for large $\kappa$ which implies \eqref{eq:EqualFields}. This finishes the proof of Theorem~\ref{thm:EqualHC3}.
\qed\\

\section{Energy of minimizers}
\label{energy}

\subsection{Basic properties}~\\
In the case where $\frac{H}{\kappa} \rightarrow \frac{1}{\mu}$, with $\Lambda_{1}=\min_{\mathsf{s}\in\Sigma}\mu_1(\alpha_{\mathsf{s}}) < \mu < \Theta_0$, superconductivity is dominated by the corners. 
The asymptotics of the ground state energy in this case is given by Theorem~\ref{thm:EnergyCorners} which we will prove in the present section.

Recall the functionals $J^{\alpha}_{\mu_1,\mu_2}$ with ground state energy $E_{\mu_1, \mu_2}^{\alpha}$ defined on angular sectors $\Gamma_{\alpha}$ by \eqref{eq:Jmu}.
We give the following proposition without proof, since it is completely analogous to the similar statements for ${\mathcal E}_{\kappa,H}$.

\begin{prop}\label{prop:AngularNL}~\\
The map $(0, \Theta_0) \times {\mathbb R}_{+} \ni (\mu_1, \mu_2) \mapsto E_{\mu_1, \mu_2}^{\alpha}$ is continuous.

Suppose that $\mu_1< \Theta_0$.
If $\mu_1 \leq \mu_1(\alpha)$, then $E_{\mu_1, \mu_2}^{\alpha}=0$ and $\psi=0$ is a minimizer.

If $\mu_1 > \mu_1(\alpha)$,
there exists a non-trivial minimizer $\psi_0$ of $J_{\mu_1, \mu_2}^{\alpha}$.
Furthermore, there exist constants $a,C>0$ such that
\begin{align}
\label{eq:DecayGamma}
\int_{\Gamma_{\alpha}} e^{2a|x|}\big( |\psi_0(x)|^2 + |(-i\nabla -{\bf F})\psi_{0}|^2\big)\,dx \leq C.
\end{align}
Finally, $\psi_0$ satisfies the uniform bound,
$$
\| \psi_0 \|_{\infty} \leq \frac{\mu_1}{\mu_2}.
$$
\end{prop}

One easily verifies the following scaling property.
\begin{prop}\label{prop:AngularNLScaling}~\\
Let $\Lambda>0$. Then the functional, 
$$
\psi \mapsto
\int_{\Gamma_{\alpha}} |(-i\nabla -\Lambda^{-2}{\bf F})\psi|^2 - \mu_1 \Lambda^{-2}|\psi |^2 + \frac{\mu_2}{2}\Lambda^{-2}|\psi|^4\,dx,
$$
defined on $\{ \psi \in L^2(\Gamma_{\alpha})\,\big|\, (-i\nabla -\Lambda^{-2}{\bf F})\psi \in L^2(\Gamma_{\alpha}) \}$ is minimized by
$\tilde \psi_0(y) = \psi_0(y/\Lambda)$, where 
$\psi_0$ is the minimizer of $J^{\alpha}_{\mu_1, \mu_2}$.

In particular,
$$
\inf_{\psi} \int_{\Gamma_{\alpha}} |(-i\nabla -\Lambda^{-2}{\bf F})\psi|^2 - \mu_1 \Lambda^{-2}|\psi |^2 + \frac{\mu_2}{2}\Lambda^{-2}|\psi|^4\,dx
=
E^{\alpha}_{\mu_1, \mu_2}.
$$
\end{prop}
By continuity of $E^{\alpha}_{\mu_1, \mu_2}$ we get the following consequence.
\begin{prop}
\label{prop:cont}
Suppose that $\frac{\kappa}{H(\kappa)} \rightarrow \mu < \Theta_0$ as $\kappa \rightarrow \infty$, and that $d_1(\kappa), d_2(\kappa) \rightarrow 1$ as $\kappa \rightarrow \infty$.
Then the ground state energy of the functional
$$
\psi \mapsto
\int_{\Gamma_{\alpha}} |(-i\nabla -\kappa H{\bf F})\psi|^2 - d_1(\kappa)\kappa^2 |\psi |^2 + d_2(\kappa)\frac{\kappa^2}{2}|\psi|^4\,dx,
$$
tends to $E^{\alpha}_{\mu, \mu}$ as $\kappa \rightarrow \infty$.
\end{prop}

\subsection{Coordinate changes}\label{change}$\,$\\
Let ${\mathsf s}\in \Sigma$. By the assumption that $\partial \Omega$ is a curvilinear domain there exists $r_{\mathsf{s}}>0$ and a local diffeomorphism $\Phi_{\mathsf{s}}$ of ${\mathbb R}^2$ such that $\Phi_{\mathsf{s}}({\mathsf s}) =0$,
$(D\Phi_{\mathsf{s}})({\mathsf s}) \in SO(2)$ and 
$\Phi_{\mathsf{s}}\big( B({\mathsf s},r_{\mathsf{s}}) \cap \Omega\big) 
= \Gamma_{\alpha_{\mathsf{s}}} \cap \Phi_{\mathsf{s}}( B({\mathsf s},r_{\mathsf{s}}))$.

Let $u, {\bf A}=(A_1,A_2) \in C^{\infty}_0( B({\mathsf s},r_{\mathsf{s}}))$ and define
$\tilde{u}(y)=u(\Phi_{\mathsf{s}}^{-1}(y))$.
Let furthermore, $\tilde{B}(y) = B(\Phi_{\mathsf{s}}^{-1}(y))$, where $B(x) = \curl {\bf A}$.
Then the quadratic form transforms as
\begin{align}
\int_{\Omega} | (-i\nabla &-{\bf A})u(x)|^2 \,dx\nonumber\\
&=
\int_{\Gamma_{\alpha_{\mathsf{s}}}} \langle (-i\nabla -{\bf \tilde A})\tilde u(y), G(y) (-i\nabla -{\bf \tilde A})\tilde u(y) \rangle\,|\det D \Phi_{\mathsf{s}}^{-1}(y)| \,dy.
\end{align}
Here
$G(y) = (D \Phi_{\mathsf{s}}) (D\Phi_{\mathsf{s}})^{T} \big|_{x = \Phi_{\mathsf{s}}(y)}$, 
and ${\bf \tilde A}=(\tilde A_1,\tilde A_2)$ satisfies
$A_1 dx_1 +  A_2 dx_2=\tilde A_1 dy_1 + \tilde A_2 dy_2$, so
\begin{align}
\label{eq:Btilde}
\partial_{y_1} \tilde A_2 - \partial_{y_2} \tilde A_1 = |\det D \Phi_{\mathsf{s}}^{-1}(y)|\; \tilde B(y).
\end{align}

\subsection{Proof of Theorem~\ref{thm:EnergyCorners}}~\\
\noindent{\bf Upper bounds}$\,$\\
We indicate here how to obtain the inequality
\begin{align}
\label{eq:UpperEnergy}
\inf_{(\psi, {\bf A})}{\mathcal E}_{\kappa, H(\kappa)}[\psi, {\bf A}] \leq \sum_{\mathsf{s}\in\Sigma} E^{\alpha_{\mathsf{s}}}_{\mu,\mu} + o(1),
\end{align}
which is the `easy' part of \eqref{eq:EnergyAsymp}. 

The inequality \eqref{eq:UpperEnergy} follows from a calculation with an explicit trial state.
The test functions will be of the form ${\bf A} = {\bf F}$ and
\begin{align*}
\psi(x) = \sum_{\mathsf{s}\in\Sigma} \psi_{\mathsf{s}}(\Phi_{\mathsf{s}}(x)),\quad\quad
\text{ with } \quad\quad
\psi_{\mathsf{s}}(y) = e^{i\kappa H \eta_{\mathsf{s}}} \psi^{\alpha_{\mathsf{s}}}_{1,1}(\sqrt{\kappa H} y) \chi(|y|).
\end{align*}
Here $\eta_{\mathsf{s}} \in C^{\infty}({\mathbb R}^2, {\mathbb R})$ is a gauge function, $\chi$ is a standard cut-off function, $\chi=1$ on a neighborhood of $0$, $\supp \chi \subset B(0,r)$, with $r = \min_{\mathsf{s}\in\Sigma} \{r_{\mathsf{s}}\}$,
and 
$\psi^{\alpha_{\mathsf{s}}}_{1,1}$ is the minimizer of $J^{\alpha_{\mathsf{s}}}_{1,1}$.
The proof of \eqref{eq:UpperEnergy} is a straight forward calculation similar to the lower bound (given below) and will be omitted.
Notice though that the decay estimates \eqref{eq:DecayGamma} for the minimizers $\psi^{\alpha_{\mathsf{s}}}_{1,1}$ imply that $\psi^{\alpha_{\mathsf{s}}}_{1,1}(\sqrt{\kappa H} y) \chi(|y|) - \psi^{\alpha_{\mathsf{s}}}_{1,1}(\sqrt{\kappa H} y)$, is exponentially small.

\noindent{\bf Lower bounds}~\\
Let $(\psi, {\bf A})$ be a minimizer of ${\mathcal E}_{\kappa, H}$.
Define $\chi_1, \chi_2 \in C^{\infty}({\mathbb R})$ to be a standard partition of unity, $\chi_1$ is non-increasing, $\chi_1^2 + \chi_2^2 =1$, $\chi_1(t) = 1$ for $t\leq 1$, $\chi_1(t) =0$ for $t\geq 2$.

For ${\mathsf{s}}\in\Sigma$, let
$$
\phi_{\mathsf{s}}(x) = \chi_1\big( \kappa^{1-\epsilon} \dist(x, {\mathsf s})\big)
$$
with $\epsilon>0$, and define $\phi_0 = \sqrt{1 - \sum_{{\mathsf{s}}\in\Sigma} \phi_{\mathsf{s}}^2}$.
Notice that when $\kappa$ is sufficiently large and ${\mathsf{s}},{\mathsf{s}}'\in\Sigma$, ${\mathsf{s}} \neq {\mathsf{s}}'$, then $\phi_{\mathsf{s}} \phi_{\mathsf{s}'} = 0$.
Therefore, using the Agmon estimates, the IMS-localization formula and the estimate $\| \psi \|_{\infty} \leq 1$, we can write,
\begin{align}
\label{eq:Divide}
{\mathcal E}_{\kappa, H}[\psi, {\bf A}] \geq \sum_{\mathsf{s}\in\Sigma} {\mathcal E}_{\kappa, H}[ \phi_{\mathsf{s}} \psi, {\bf A}] +
{\mathcal O}(\kappa^{-\infty}).
\end{align}
By the Sobolev imbedding $W^{1,2}(\Omega) \rightarrow L^4(\Omega)$, Lemma~\ref{lem:H1ControlsL2} combined with \eqref{eq:BasicCurl}, and the Agmon estimate we get
\begin{align}
(\kappa H)^2 \| {\bf A} - {\bf F} \|_{4}^2 &
\leq C (\kappa H)^2 \| {\bf A} - {\bf F} \|_{W^{1,2}(\Omega)}^2 \nonumber\\
&\leq C' (\kappa H)^2 \| \curl {\bf A} - 1 \|_{2}^2 \nonumber\\
&\leq C'' \kappa^2 \| \psi \|_2^2 \leq C'''.
\end{align}
Thus we can estimate
\begin{align}
&\int_{\Omega} \big| (-i\nabla - \kappa H {\bf A}) (\phi_{\mathsf{s}} \psi) \big|^2\,dx\nonumber\\
&\geq
(1- \kappa^{-1/2}) \int_{\Omega} \big| (-i\nabla - \kappa H {\bf F}) (\phi_{\mathsf{s}} \psi) \big|^2\,dx 
- \kappa^{1/2} (\kappa H)^2 \| {\bf A} - {\bf F} \|_{4}^2 \| \phi_{\mathsf{s}} \psi \|_{4}^2\nonumber\\
&\geq
(1- \kappa^{-1/2}) \int_{\Omega} \big| (-i\nabla - \kappa H {\bf F}) (\phi_{\mathsf{s}} \psi) \big|^2\,dx 
-C \kappa^{-1/2},
\end{align}
where we used the inequality
$$
\| \phi_{\mathsf{s}} \psi \|_{4}^2\leq
\sqrt{ \int_{\Omega} |\psi|^2\,dx}\leq
\sqrt{ C \int_{\{\dist(x, \Sigma)\leq M \kappa^{-1}\}} 1 \,dx} \leq C' \kappa^{-1}.
$$
Now consider the change of coordinates $\Phi_{\mathsf{s}}$ from subsection~\ref{change}. 
For sufficiently large values of $\kappa$ we have
$\supp \phi_{\mathsf{s}} \subset B({\mathsf s}, r_{\mathsf{s}})$.
Define 
$$
\tilde \psi_{\mathsf{s}} = (\phi_{\mathsf{s}} \psi) \circ \Phi_{\mathsf{s}}^{-1}.
$$
Since $| \det D \Phi_{\mathsf{s}}(0) | = 1$, we get by Taylor's formula that
$$
\big| \, | \det D \Phi_{\mathsf{s}} | - 1 \big| \leq C \kappa^{-1 + \epsilon},
$$
on $\supp \tilde \psi_{\mathsf{s}}$.

Consider the transformed magnetic field as in \eqref{eq:Btilde}.
We define
\begin{align}
\tilde \beta(y) :=  |\det D \Phi_{\mathsf{s}}^{-1}(y)| \tilde B(y)
= |\det D \Phi_{\mathsf{s}}^{-1}(y)| = 1 + {\mathcal O}(\kappa^{-1 + \epsilon}),
\end{align}
on $\supp \tilde \psi_{\mathsf{s}}$.
We look for ${\bf \tilde A}=(\tilde A_1,\tilde A_2)$ such that
$\partial_{y_1} \tilde A_2 - \partial_{y_2} \tilde A_1 = \tilde \beta(y)$.

One choice of a solution is
$$
{\bf \tilde A} = \left(-y_2/2, \int_0^{y_1} [\tilde \beta(y_1', y_2) -1/2]\,dy_1'\right).
$$
With this choice
$$
\| {\bf \tilde A}- {\bf F} \|_{L^{\infty}(B(0, C \kappa^{-1 + \epsilon}))} \leq 
C' \kappa^{-2 + 2\epsilon}.
$$
Thus
\begin{align}
(\kappa H)^2 \int |{\bf \tilde A} - {\bf F}|^2 \,|\tilde \psi_{\mathsf{s}}|^2\,dy
\leq
C (\kappa H)^2 \kappa^{-4 + 4\epsilon} \int |\tilde \psi_{\mathsf{s}}|^2\,dy
\leq C' \kappa^{-2 + 6 \epsilon}.
\end{align}
Therefore, for some $\eta \in C^{\infty}(\Omega, {\mathbb R})$ we find
\begin{align}
\label{eq:88}
&\int_{\Omega} \big| (-i\nabla - \kappa H {\bf F}) (\phi_{\mathsf{s}} \psi) \big|^2\,dx\nonumber\\
&=
\int_{\Gamma_{\alpha_{\mathsf{s}}}} \big\langle (-i\nabla - \kappa H {\bf \tilde A})( e^{i \kappa H \eta} \tilde \psi_{\mathsf{s}}),
G(y) (-i\nabla - \kappa H {\bf \tilde A})( e^{i \kappa H \eta} \tilde \psi_{\mathsf{s}}) \big\rangle | \det D \Phi_{\mathsf{s}}^{-1} | \,dy\nonumber\\
&\geq
(1 - C \kappa^{-1+\epsilon})
\int_{\Gamma_{\alpha_{\mathsf{s}}}} \big | (-i\nabla - \kappa H {\bf \tilde A})( e^{i \kappa H \eta} \tilde \psi_{\mathsf{s}})
\big|^2 \,dy\nonumber\\
&\geq
(1 - C \kappa^{-1+\epsilon})
\Big\{ (1- \kappa^{-1+3\epsilon})
\int_{\Gamma_{\alpha_{\mathsf{s}}}} \big | (-i\nabla - \kappa H {\bf F})( e^{i \kappa H \eta} \tilde \psi_{\mathsf{s}})
\big|^2 \,dy\nonumber\\
&\quad\quad\quad\quad\quad\quad\quad\quad\quad\quad\quad\quad-
\kappa^{1-3\epsilon} (\kappa H)^2 \int_{\Gamma_{\alpha_{\mathsf{s}}}} | {\bf \tilde A} - {\bf F}|^2 |\tilde \psi_{\mathsf{s}}|^2\,dy
\Big\}\nonumber\\
&\geq
(1 - 2C \kappa^{-1+3\epsilon})
\int_{\Gamma_{\alpha_{\mathsf{s}}}} \big | (-i\nabla - \kappa H {\bf F})( e^{i \kappa H \eta} \tilde \psi_{\mathsf{s}})
\big|^2 \,dy + {\mathcal O}(\kappa^{-1 + 3 \epsilon}).
\end{align}
By \eqref{eq:88} we find
\begin{align}
\label{eq:Endelig}
{\mathcal E}_{\kappa, H}&[ \phi_{\mathsf{s}} \psi, {\bf A}] \nonumber\\
&\geq
(1- C_1 \kappa^{-1+3\epsilon})
\int_{\Gamma_{\alpha_{\mathsf{s}}}} \big\{ \big | (-i\nabla - \kappa H {\bf F})( e^{i \kappa H \eta} \tilde \psi_{\mathsf{s}})
\big|^2\nonumber\\
&\quad\quad\quad\quad\quad\quad\quad\quad\quad\quad
- (1+C_2 \kappa^{-1+3\epsilon}) \frac{\kappa^2}{2} |  e^{i \kappa H \eta} \tilde \psi_{\mathsf{s}}|^2
+ \kappa^4| e^{i \kappa H \eta} \tilde \psi_{\mathsf{s}}|^4\big\}\,dy\nonumber\\
&\quad\quad\quad\quad+ {\mathcal O}(\kappa^{-1 + 3 \epsilon}).
\end{align}
We choose $0<\epsilon<1/3$ arbitrary.
Using Proposition~\ref{prop:cont} and combing \eqref{eq:Divide} and \eqref{eq:Endelig} we find the lower bound inherent in \eqref{eq:UpperEnergy}, i.e.
$$
{\mathcal E}_{\kappa, H(\kappa)}[\psi, {\bf A}] \geq 
\sum_{\mathsf{s}\in\Sigma} E^{\alpha_{\mathsf{s}}}_{\mu,\mu} + o(1).
$$
This finishes the proof of Theorem~\ref{thm:EnergyCorners}.
\qed

\medskip\par
\noindent{\bf Acknowledgements}\\
It is a pleasure to acknowledge discussions on this and related subjects with X.~Pan and R.~Frank. Furthermore, without the discussions with and encouragement from B.~Helffer this work would never have been carried through.\\
Both authors were supported by the ESF
Scientific Programme in Spectral Theory and Partial Differential
Equations (SPECT).
The second author is supported by a Skou-stipend from the Danish Research Council and has also benefitted from support from the European Research Network
`Postdoctoral Training Program in Mathematical Analysis of Large
Quantum Systems' with contract number HPRN-CT-2002-00277. 
Furthermore, the second author wants to thank CIMAT in Guanajuato, Mexico for hospitality.

\end{document}